\documentclass[11pt,amsfonts]{amsart}

\usepackage{fullpage}

\newtheorem{theorem}{Theorem}[section]

\newtheorem{corollary}[theorem]{Corollary}

\newtheorem{lemma}[theorem]{Lemma}

\newtheorem{proposition}[theorem]{Proposition}

\newtheorem{remark}[theorem]{Remark}

\newtheorem{example}[theorem]{Example}

\def\endproof{\qed \medskip}
\def\blacksquare{\hbox to .60em{\vrule width .60em height .60em}}

\newcommand{\Met}{\operatorname{Met}}
\newcommand{\Diff}{\operatorname{Diff}}
\newcommand{\Ric}{\operatorname{Ric}}

\newcounter{mnotecount}[section]

\begin{document}

\title[]{Unique continuation results for Ricci curvature and 
Applications}

\author[]{Michael T. Anderson and Marc Herzlich}

\thanks{The first author is partially supported by NSF Grant DMS
 0604735; the second author 
is partially supported by ANR project GeomEinstein 06-BLAN-0154.  }

\abstract{Unique continuation results are proved for metrics with
 prescribed 
Ricci curvature in the setting of bounded metrics on compact manifolds
 with 
boundary, and in the setting of complete conformally compact metrics on
 such 
manifolds. Related to this issue, an isometry extension property is
 proved: 
continuous groups of isometries at conformal infinity extend into the
 bulk of 
any complete conformally compact Einstein metric. Relations of this
 property 
with the invariance of the Gauss-Codazzi constraint equations under
 deformations 
are also discussed.}
\endabstract

\maketitle

\setcounter{section}{0}

\section{Introduction.}
\setcounter{equation}{0}

 In this paper, we study certain issues related to the boundary 
behavior of metrics with prescribed Ricci curvature. Let $M$ be a 
compact $(n+1)$-dimensional manifold with compact non-empty boundary 
$\partial M$. We consider two possible classes of Riemannian metrics 
$g$ on $M$. First, $g$ may extend smoothly to a Riemannian metric on 
the closure $\bar M = M \cup  \partial M$, thus inducing a Riemannian 
metric $\gamma  = g|_{\partial M}$ on $\partial M$. Second, $g$ may be 
a complete metric on $M$, so that $\partial M$ is ``at infinity''. In 
this case, we assume that $g$ is conformally compact, i.e.~there exists
a defining function $\rho$ for $\partial M$ in $M$ such that the 
conformally equivalent metric
\begin{equation} \label{e1.1}
\widetilde g = \rho^{2}g 
\end{equation}
extends at least $C^{2}$ to $\partial M$. The defining function $\rho$ is 
unique only up to multiplication by positive functions; hence only the 
conformal class $[\gamma]$ of the associated boundary metric $\gamma  =
\bar g|_{\partial M}$ is determined by $(M, g)$. 

 The issue of boundary regularity of Riemannian metrics $g$ with 
controlled Ricci curvature has been addressed recently in several 
papers. Thus, [4] proves boundary regularity for bounded metrics $g$ 
on $M$ with controlled Ricci curvature, assuming control on the 
boundary metric $\gamma$ and the mean curvature of $\partial M$ in $M$.
In [16], boundary regularity is proved for conformally compact Einstein
metrics with smooth conformal infinity; this was previously proved by 
different methods in dimension 4 in [3], cf.~also [5].

\medskip

 One purpose of this paper is to prove a unique 
continuation property at the boundary $\partial M$ for bounded metrics
 or for 
conformally compact metrics. We first state a version of the result for
 Einstein 
metrics on bounded domains.
\begin{theorem} \label{t 1.1.}
  Let $(M, g)$ be a $C^{3,\alpha}$ metric on a compact manifold with
 boundary $M$, with induced metric $\gamma  = g|_{\partial M}$, 
and let $A$ be the $2^{\rm nd}$ fundamental form of $\partial M$ 
in $M$. Suppose the Ricci curvature $\Ric_{g}$ satisfies
\begin{equation} \label{e1.2}
\Ric_{g} = \lambda g, 
\end{equation}
where $\lambda$ is a fixed constant.

 Then $(M, g)$ is uniquely determined up to local isometry and inclusion, 
by the Cauchy data $(\gamma, A)$ on an arbitrary open set $U$ of $\partial M$. 
\end{theorem}

  Thus, if $(M_{1}, g_{1})$ and $(M_{2}, g_{2})$ are a pair of Einstein metrics 
as above, whose Cauchy data $(\gamma, A)$ agree on an open set $U$ common to 
both $\partial M_{1}$ and $\partial M_{2}$, then after passing to suitable 
covering spaces $\bar M_{i}$, either there exist isometric embeddings 
$\bar M_{1} \subset \bar M_{2}$ or $\bar M_{2} \subset \bar M_{1}$ or there 
exists an Einstein metric $(\bar M_{3}, g_{3})$ and isometric embeddings 
$(\bar M_{i}, g_{i}) \subset (\bar M_{3}, g_{3})$. Similar results hold for 
metrics which satisfy other covariant equations involving the metric to 
$2^{\rm nd}$ order, for example the Einstein equations coupled to other 
fields; see Proposition 3.7. 

\medskip

  For conformally compact metrics, the $2^{\rm nd}$ fundamental form
 $A$ of the 
compactified metric $\bar g$ in (1.1) is umbilic, and completely
 determined by 
the defining function $\rho$. In fact, for conformally compact Einstein
 metrics, 
the higher order Lie derivatives ${\mathcal L}_{N}^{(k)}\bar g$ at
 $\partial M$, where 
$N$ is the unit vector in the direction $\bar \nabla \rho$, are
 determined by the 
conformal infinity $[\gamma]$ and $\rho$ up to order $k < n$. Supposing
 $\rho$ is a 
geodesic defining function, so that $||\bar \nabla \rho|| = 1$, let 
\begin{equation}\label{e1.3}
g_{(n)} = \tfrac{1}{n!}{\mathcal L}_{N}^{(n)}\bar g.
\end{equation}
More precisely, $g_{(n)}$ is the $n^{\rm th}$ term in the
 Fefferman-Graham expansion 
of the metric $g$; this is given by (1.3) when $n$ is odd, and in a
 similar way when $n$ 
is even, cf. [18] and \S 4 below. The term $g_{(n)}$ is the natural
 analogue of $A$ for 
conformally compact Einstein metrics.

\begin{theorem} \label{t 1.2.}
  Let $g$ be a $C^{2}$ conformally compact Einstein metric on a compact
 manifold $M$ with 
$C^{\infty}$ smooth conformal infinity $[\gamma]$, normalized so that
\begin{equation} \label{e1.4}
\Ric_{g} = -ng, 
\end{equation}
 Then the Cauchy data $(\gamma, g_{(n)})$ restricted to any open set
 $U$ of 
$\partial M$ uniquely determine $(M, g)$ up to local isometry and
 determine 
$(\gamma, g_{(n)})$ globally on $\partial M$.
\end{theorem}

  The recent boundary regularity result of Chru\'sciel et al.,~[16],
 implies that $(M, g)$ 
is $C^{\infty}$ polyhomogeneous conformally compact, so that the
 hypotheses of Theorem 1.2 
imply the term $g_{(n)}$ is well-defined on $\partial M$. A more
 general version of Theorem 1.2, without the smoothness assumption 
on $[\gamma]$, is proved in \S 4, cf.~Theorem 4.1. For conformally 
compact metrics coupled to other fields, see Remark 4.5. 

\medskip

  Of course neither Theorem 1.1 or 1.2 hold when just the boundary
 metric 
$\gamma$ on $U \subset \partial M$ is fixed. For example, in the
 context of 
Theorem 1.2, by [20] and [16], given any $C^{\infty}$ smooth boundary
 metric $\gamma$ 
sufficiently close to the round metric on $S^{n}$, there is a smooth
 (in the polyhomogeneous 
sense) conformally compact Einstein metric on the $(n+1)$-ball
 $B^{n+1}$, close to the 
Poincar\'e metric. Hence, the behavior of $\gamma$ in $U$ is
 independent of its 
behavior on the complement of $U$ in $\partial M$. 

  Theorems 1.1 and 1.2 have been phrased in the context of ``global''
 Einstein 
metrics, defined on compact manifolds with compact boundary. However,
 the proofs 
are local, and these results hold for metrics defined on an open
 manifold with 
boundary. From this perspective, the data $(\gamma, A)$ or $(\gamma,
 g_{(n)})$ 
on $U$ determine whether Einstein metric $g$ has a global extension to
 an 
Einstein metric on a compact manifold with boundary, (or conformally
 compact Einstein 
metric), and how smooth that extension is at the global boundary. 

\medskip

  A second purpose of the paper is to prove the following isometry
 extension result which 
is at least conceptually closely related to Theorem 1.2. However, while
 Theorem 1.2 is 
valid locally, this result depends crucially on global properties. 

\begin{theorem} \label{t1.3}
Let $g$ be a $C^{2}$ conformally compact Einstein metric on a compact
 manifold $M$ with 
$C^{\infty}$ boundary metric $(\partial M, \gamma)$, and suppose
\begin{equation}\label{e1.5}
\pi_{1}(M, \partial M) = 0.
\end{equation}
Then any connected group of isometries of $(\partial M, \gamma)$
 extends to an action 
by isometries on $(M, g)$. 
\end{theorem}

  The condition (1.5) is equivalent to the statement that $\partial M$
 is connected 
and  the inclusion map $\iota: \partial M \rightarrow M$ induces a
 surjection
$\pi_{1}(\partial M) \rightarrow \pi_{1}(M) \rightarrow 0$.

  Rather surprisingly, this result is closely related to the equations
 at conformal 
infinity induced by the Gauss-Codazzi equations on hypersurfaces
 tending to $\partial M$. 
It turns out that isometry extension from the boundary at least into a
 thickening of 
the boundary is equivalent to the requirement that the Gauss-Codazzi
 equations induced 
at $\partial M$ are preserved under arbitrary deformations of the
 boundary metric. This 
is discussed in detail in \S 5, see e.g.~Proposition 5.4. We note that
 this result does 
not hold for complete, asymptotically (locally) flat Einstein metrics,
 cf. Remark 5.8. 

  A simple consequence of Theorem 1.3 is the following uniqueness
 result: 
\begin{corollary} \label{c1.4}
A $C^{2}$ conformally compact Einstein metric with conformal infinity
 given by the 
class of the round metric $g_{+1}$ on the sphere $S^{n}$ is necessarily
 isometric 
to the Poincar\'e metric on the ball $B^{n+1}$. 
\end{corollary}

  Results similar to Theorem 1.3 and Corollary 1.4 have previously been
 proved in a number 
of different special cases by several authors, see for example [7],
 [9], [31], [33]; 
the proofs in all these cases are very different from the proof given
 here. 

\medskip

  It is well-known that unique continuation does not hold for large
 classes of elliptic 
systems of PDE's, even for general small perturbations of systems which
 are diagonal at 
leading order; see for instance [23] and references therein for a
 discussion related to 
geometric PDEs. The proofs of Theorems 1.1 and 1.2 rely on unique
 continuation results 
of Calder\'on [13], [14] and Mazzeo [27] respectively, based on
 Carleman estimates. The 
main difficulty in reducing the proofs to these results is the
 diffeomorphism covariance 
of the Einstein equations and, more importantly, that of the
 ``abstract'' Cauchy data 
$(\gamma, A)$ or $(\gamma, g_{(n)})$ at $\partial M$. The unique
 continuation theorem 
of Mazzeo requires a diagonal (i.e.~uncoupled) Laplace-type system of
 equations, at 
leading (second) order. The unique continuation result of Calder\'on is
 more general, 
but again requires strong restrictions on the structure of the leading
 order symbol 
of the operator. For emphasis and clarity, these issues are discussed
 in more detail 
in \S 2. The proofs of Theorems 1.1, 1.2 and 1.3 are then given in \S
 3, \S 4 and \S 5 
respectively. 

\medskip

Very recently, while the writing on this paper was being completed, 
O. Biquard [12] has given a different proof of Theorem 1.2, which avoids 
some of the gauge issues discussed above. However, his method apparently 
requires $C^{\infty}$ smoothness of the boundary data, which limits the 
applicability of this result; for instance the applications in [5] 
or [6] require finite or low differentiability of the boundary data. 

\medskip

We would like to thank Michael Taylor for interesting discussions 
on geodesic-harmonic coordinates, Piotr Chru\'sciel and Erwann Delay for
interesting discussions concerning Theorem 1.3, and Olivier Biquard for
informing us of his independent work on unique continuation.

\section{Local Coordinates and Cauchy Data}
\setcounter{equation}{0}

In this section, we discuss in more detail the remarks in the Introduction 
on classes of local coordinate systems, and their relation with Cauchy
data on the boundary $\partial M$. 

Thus, consider for example solutions to the system
\begin{equation} \label{e2.1}
\Ric_{g} = 0, 
\end{equation}
defined near the boundary $\partial M$ of an $(n+1)$-dimensional manifold $M$. 
Since the Ricci curvature involves two derivatives of the metric, 
Cauchy data at $\partial M$ consist of the boundary metric $\gamma $ 
and its first derivative, invariantly represented by the $2^{\rm nd}$ 
fundamental form $A$ of $\partial M$ in $M$. Thus, we assume $(\gamma , A)$ 
are prescribed at $\partial M$, (subject to the Gauss and Gauss-Codazzi
equations), and call $(\gamma, A)$ abstract Cauchy data. Observe that 
the abstract Cauchy data are invariant under diffeomorphisms of $M$ equal 
to the identity at $\partial M$. 

\medskip

The metric $g$ determines the geodesic defining function 
$$t(x) = dist_{g}(x, \partial M).$$
The function $t$ depends of course on $g$; however, given any 
other smooth metric $g'$, there is a diffeomorphism $F$ of a 
neighborhood of $\partial M$, equal to the identity on $\partial M$, 
such that $t'(x) = dist_{F^{*}g'}(x, \partial M)$ satisfies $t' = t$. 
As noted above, this normalization does not change the abstract Cauchy
data $(\gamma, A)$ and preserves the isometry class of the metric. 

Let $\{y^{\alpha}\}$, $0 \leq \alpha \leq n$, be any local coordinates 
on a domain $\Omega $ in $M$ containing a domain $U$ in $\partial M$. 
We assume that $\{y^{i}\}$ for $1 \leq i \leq n$ form local coordinates 
for $\partial M$ when $y^{0} = 0$, so that $\partial / \partial y^{0}$ is 
transverse to $\partial M$. Throughout the paper, Greek indices $\alpha$, 
$\beta$ run from $0$ to $n$, while Latin indices $i$, $j$ run from $1$ 
to $n$. If $g_{\alpha\beta}$ are the components of $g$ in these coordinates, 
then the abstract Cauchy problem associated to (2.1) in the local coordinates 
$\{y^{\alpha}\}$ is the system 
\begin{equation} \label{e2.2}
(\Ric_{g})_{\alpha\beta} = 0, \ \ {\rm with} \ \  
g_{ij}|_{U} = \gamma_{ij}, \ \tfrac{1}{2}({\mathcal L}_{\nabla
 t}g)_{ij}|_{U} = 
a_{ij},
\end{equation}
where $\gamma_{ij}$ and $a_{ij}$ are given on $U$, (subject to the constraints 
of the Gauss and Gauss-Codazzi equations). Here one immediately sees a 
problem, in that (2.2) on $U \subset \partial M$ involves only the
tangential part $g_{ij}$ of the metric (at 0 order), and not the full metric 
$g_{\alpha\beta}$ at $U$. The normal $g_{00}$ and mixed $g_{0i}$ components 
of the metric are not prescribed at $U$. As seen below, these components 
are gauge-dependent; they cannot be prescribed ``abstractly'', independent 
of coordinates, as is the case with $\gamma$ and $A$.
In other words, if (2.1) is expressed in local coordinates $\{y^{\alpha}\}$ 
as above, then a well-defined Cauchy or unique continuation problem has
the form 
\begin{equation} \label{e2.3}
(\Ric_{g})_{\alpha\beta} = 0, \ \ {\rm with} \ \  
g_{\alpha\beta} = \gamma_{\alpha\beta}, \ 
\tfrac{1}{2}\partial_{t}g_{\alpha\beta} = a_{\alpha\beta}, \ {\rm on} \
  U \subset 
\partial M , 
\end{equation}
where $\Omega$ is an open set in $({\mathbb R}^{n+1})^{+}$ with
 $\partial \Omega = U$ an 
open set in $\partial ({\mathbb R}^{n+1})^{+} ={\mathbb R}^{n}$.
Formally, (2.3) is a determined system, while (2.2) is underdetermined. 

\medskip

Let $g_{0}$ and $g_{1}$ be two solutions to (2.1), with the same 
Cauchy data $(\gamma , A)$, and with geodesic defining functions 
$t_{0}$, $t_{1}$. Changing the metric $g_{1}$ by a diffeomorphism 
if necessary, one may assume that $t_{0} = t_{1}$. One may then write
the metrics with respect to a Gaussian or geodesic boundary coordinate 
system $(t, y^{i})$ as 
\begin{equation} \label{e2.4}
g_{k} = dt^{2} + (g_{k})_{t}, 
\end{equation}
where $(g_{k})_{t}$ is a curve of metrics on $\partial M$ and $k = 0, 1$. 
Here $y_{i}$ are coordinates on $\partial M$ which are extended into $M$ 
to be invariant under the flow of the vector field $\nabla t$. The metric
$(g_{k})_{t}$ is the metric induced on $S(t)$ and pulled back to 
$\partial M$ by the flow of $\nabla t$. One has $(g_{k})_{0} = \gamma$ and
$\frac{1}{2}\frac{d}{dt}(g_{k})_{t}|_{t=0} = A$. 
Since $g_{0\alpha} = \delta_{0\alpha}$ in these coordinates, 
$\nabla t = \partial_{t}$, 
and hence the local coordinates are the same for both metrics, (or at least 
may be chosen to be the same). Thus, geodesic boundary coordinates are natural
from the point of view of the Cauchy or unique continuation problem, since 
in such local coordinates the system (2.2), together with the prescription
$g_{0\alpha} = \delta_{0\alpha}$, is equivalent to the system (2.3).
However, the Ricci curvature is not elliptic or diagonal to leading order 
in these coordinates. The expression of the Ricci curvature in such
coordinates does not satisfy the hypotheses of Calder\'on's theorem [14], 
and it appears to be difficult to establish unique continuation of solutions 
in these coordinates by working directly on the equations on the metric
(see, however, [12] for another approach).

Next suppose that $\{x^{\alpha}\}$ are boundary harmonic coordinates,
defined as follows. For $1 \leq  i \leq n$, let $\hat x^{i}$ be local harmonic 
coordinates on a domain $U$ in $(\partial M, \gamma)$. Extend $\hat x^{i}$ into 
$M$ to be harmonic functions in $(\Omega, g)$, $\Omega \subset  M$,
with Dirichlet boundary data; thus
\begin{equation} \label{e2.5}
\Delta_{g}x^{i} = 0, \ \ x^{i}|_{U} = \hat x^{i}.
\end{equation}
Let $x^{0}$ be a harmonic function on $\Omega$ with 0 boundary data, so
 that
\begin{equation} \label{e2.6}
\Delta_{g}x^{0} = 0, \ \ x^{0}|_{U} = 0.
\end{equation}
Then the collection $\{x^{\alpha}\}$, $0 \leq \alpha \leq  n$, form a 
local harmonic coordinate chart on a domain $\Omega \subset (M, g)$. In
such coordinates, one has
\begin{equation} \label{e2.7}
(\Ric_{g})_{\alpha\beta} =
 -\tfrac{1}{2}g^{\mu\nu}\partial_{\mu}\partial_{\nu}g_{\alpha\beta} 
+ Q_{\alpha\beta}(g, \partial g), 
\end{equation}
where $Q(g, \partial g)$ depends only on $g$ and its first derivatives.
This is an elliptic operator, diagonal at leading order, and satisfies the
hypotheses of Calder\'on's theorem. However, in general, the local Cauchy problem
(2.3) is not well-defined in these coordinates; if $g_{0}$ and $g_{1}$ are two
solutions of (2.1), each with corresponding local boundary harmonic 
coordinates, then the components 
$(g_{0})_{0\alpha}$ and $(g_{1})_{0\alpha}$ in general will differ at
$U \subset 
\partial M$. This is of course closely related to the fact that there
 are many possible 
choices of harmonic functions $x^{\alpha}$ satisfying (2.5) and (2.6),
 and to the fact 
that the behavior of harmonic functions depends on global properties of
 $(\Omega, g)$. 
In any case, it is not known how to set up a well-defined Cauchy
 problem in these 
coordinates for which one can apply standard unique continuation
 results. 

  Consider then geodesic-harmonic coordinates ``intermediate'' between
 geodesic 
boundary and boundary harmonic coordinates. Thus, let $t$ be the
 geodesic distance 
to $\partial M$ as above. Choose local harmonic coordinates $\hat
 x^{i}$ on 
$\partial M$ as before and extend them into $M$ to be harmonic on the
 level sets 
$S(t)$ of $t$, i.e.~locally on $S(t)$, 
\begin{equation} \label{e2.8}
\Delta_{U(t)}x^{i} = 0, \ \ x^{i}|_{\partial U(t)} = \hat
 x^{i}|_{\partial U(t)}; 
\end{equation}
here the boundary value $\hat x^{i}$ is the extension of $\hat x^{i}$
 on $U$ 
into $M$ which is invariant under the flow $\phi_{t}$ of $\nabla t$, 
and $U(t) = \phi_{t}(U) \subset S(t)$. The functions $(t, x^{i})$ form
 a coordinate 
system in a neighborhood $\Omega$ in $M$ with $\Omega\cap\partial M =
 U$. 

  It is not difficult to prove that geodesic-harmonic coordinates
 preserve the 
Cauchy data, in the sense that the data (2.2) in such coordinates imply
 the data 
(2.3). However, the Ricci curvature is not an elliptic operator in the
 metric 
in these coordinates, nor is it diagonal at leading order; the main
 reason is 
that the mean curvature of the level sets $S(t)$ is not apriori
 controlled. 
So again, it remains an open question whether unique continuation can
 be proved 
in these coordinates. 

  Having listed these attempts which appear to fail, a natural choice
 of 
coordinates which do satisfy the necessary requirements are
 $H$-harmonic 
coordinates $(\tau , x^{i})$, whose $\tau$-level surfaces
 $S_{\tau}$ are of prescribed mean curvature $H$ and with $x^{i}$ 
harmonic on $S_{\tau}$. These coordinates were introduced by 
Andersson-Moncrief [8] to prove a well-posedness result for the 
Cauchy problem for the Einstein equations in general relativity, 
and, as shown in [8], have a number of advantageous properties. 
Thus, adapting some of the arguments of [8], we show in \S 3 that 
the Einstein equations (1.2) are effectively elliptic in such
 coordinates, and 
such coordinates preserve the Cauchy data in the sense above,
 (i.e.~(2.2) implies (2.3)). It will then be shown that unique 
continuation holds in such coordinates, via application of the 
Calder\'on theorem.

%
%
%
%

\section{Proof of Theorem 1.1}
\setcounter{equation}{0}

Theorem 1.1 follows from a purely local result, which we formulate as follows. 
Let $C$ be a domain diffeomorphic to a cylinder $I\times B^{n} \subset 
{\mathbb R}^{n+1}$, with $U = \{0\}\times B^n$, diffeomorphic to a ball 
in ${\mathbb R}^{n}$. Let $U = \partial_{0}C$ be the horizontal boundary 
and $\partial C = I\times S^{n-1}$ be the vertical boundary. 

Let $g$ be a Riemannian metric on $C$ which is $C^{k-1,\alpha}$ up to the 
boundary of $C$ in the given standard coordinate 
system $\{y^{\alpha}\} = \{y^0,y^i\}$ with $y^{0} = 0$ on $U$ and $k\geq 2$. 
Without loss of generality, we assume that $C$ is chosen sufficiently 
small so that $g$ is close to the Euclidean metric $\delta$ in the
 $C^{k-1,\alpha}$ topology. For simplicity, we shall rescale $C$ and 
the coordinates $\{y^{\alpha}\}$ if necessary so that $(C, g)$ is 
$C^{k-1,\alpha}$ close to the standard cylinder $((I\times B^{n}(1), 
B^{n}(1))) \subset ({\mathbb R}^{n+1},{\mathbb R}^{n})$, $I = [0,1]$. 

We will prove the following local version of Theorem 1.1. 

\begin{theorem} \label{t 3.1.}
Let $g_{0}$, $g_{1}$ be two $C^{k-1,\alpha}$ metrics as above on $C$, 
$k\geq 4$, satisfying
\begin{equation} \label{e3.1}
\Ric_{g_{i}} = \lambda g_{i}, \quad i = 0,1
\end{equation}
for some fixed constant $\lambda$. Suppose $g_{0}$ and 
$g_{1}$ have the same abstract Cauchy data on $U$ in the sense of \S 2,
so that $\gamma_{0} = \gamma_{1}$ and $A_{0} = A_{1}$. 

Then $(C, g_{0})$ is isometric to $(C, g_{1})$, by an isometry 
equal to the identity on $U$. In particular, Theorem 1.1
holds. 
\end{theorem}

The proof of Theorem 3.1 will proceed in several steps, organized 
around several Lemmas. We first work with a fixed 
metric $g$ on $C$ as above. Let $N$ be the inward unit normal to $U$ in 
$C$ and let $A = \nabla N$ be the corresponding second fundamental 
form, with mean curvature $H = tr_g A$ on $U$. By the initial assumptions
above, $A$ and $H$ are close to $0$ in $C^{k-2,\alpha}$; more precisely,
one may assume that 
$$||A||_{C^{k-2,\alpha}} = O(\varepsilon)$$
with $\varepsilon$ positive but as small as needed, by a further rescaling 
of the coordinates (this will play an important role at various places
below). Note moreover that the rescaling process turns the Einstein constant 
$\lambda$ into $\varepsilon\lambda$. Abusing notation here, we denote $y^{0} = t$ 
and without loss of generality assume that the coordinates $y^{i}$ are harmonic 
on $U$. 

\medskip

To begin, we construct certain systems of $H$-harmonic coordinates discussed 
at the end of \S 2. Let $\phi: C \rightarrow C$ be a diffeomorphism of the 
cylinder $C$, (in other words a change of coordinates), so that 
$y^{\alpha} = \phi^{\alpha}(x^{\beta})$, where $x^{\beta}$ is another 
coordinate system for $C$. As above, we write $x^{\alpha} = (\tau, x^{i})$ 
and assume that $\phi$ is close to the identity map. The level surfaces 
$S_{\tau} = \{\tau\}\times B^{n}$ are mapped under $\phi$ to a foliation 
$\Sigma_{\tau}$ of $C$, with each leaf given by the graph of the function 
$\phi_{\tau}$ over $B^{n}$. We assume $\phi_{0} = id$, so that $\phi = id$ 
on $U$. Let $f:\partial C \rightarrow \partial C$ be the induced diffeomorphism 
on the boundary $\partial C$. 

\begin{lemma} \label{l 3.2.}
Let $k\geq 2$. Given a $C^{k,\alpha}$ mapping $f$ on $\partial C$ as above, 
close to the identity in $C^{k,\alpha}$, and a metric $g$ close to the 
Euclidean metric $\delta$ in $C^{k-1,\alpha}$ on $C$, there exists a 
unique $\phi \in \Diff^{k,\alpha}(C)$ such that, with respect to the 
pull-back metric $\phi^{*}(g)$, 
\begin{equation} \label{e3.2}
H^{\phi^{*}(g)}(S_{\tau}) = H^{\phi^{*}(g)}(S_{0}), \ \ {\rm and} \ \ 
\Delta^{\phi^{*}(g)}_{S_{\tau}}x^{i} = 0,
\end{equation}
with the property that $\phi_{|\partial C} = f$. Thus, the leaves 
$\tau = const$ have mean curvature independent of $\tau$, in the 
$x^{\alpha}$-coordinates, and the coordinate functions $x^{i}$ are 
harmonic on each $S_{\tau}$. 
\end{lemma}

{\bf Proof:} Let 
$$ \mathcal{H} : \Met^{k-1,\alpha}(C)\times \Diff^{k,\alpha}_0(C) 
\longrightarrow C^{k-2,\alpha}(C)\times \prod_{1}^{n}C^{k-2,\alpha}(C) 
\times \Diff^{k,\alpha}_0(\partial C) $$
$$ \mathcal{H}(g,F) = (H^{F^{*}(g)}(S_{\tau}) - H^{F^{*}(g)}(S_{0}), 
\Delta^{F^{*}(g)}_{S_{\tau}}x^{i}, F_{|\partial C}),$$
where $\Diff_{0}^{k,\alpha}(C)$ is the space of $C^{k,\alpha}$ diffeomorphisms 
on the cylinder equal to the identity on $C_{0} = \{0\}\times B^{n}$. 
The map $\mathcal{H}$ is clearly a smooth map of Banach spaces, and its
 linearization at $(\delta, id)$ in the second variable is
$$ L(v) = (\Delta_{\delta}v^{0}, \Delta_{\delta}v^{i}, v_{|\partial C}), $$
where $\Delta_{\delta}$ is the Laplacian with respect to the flat metric $\delta$ 
on $S_{\tau}$. The operator $L$ is clearly an isomorphism, and by the implicit 
function theorem in Banach spaces, it follows that there is a smooth map
$$\Phi : \mathcal{U}\times\mathcal{V} \subset \Met^{k-1,\alpha}(C)\times
 \Diff^{k,\alpha}_0(\partial C) 
\longrightarrow \Diff^{k,\alpha}_0(C),$$
$$\Phi(g,f) = \phi^g(f)$$
from a neigbourhood of the Euclidean metric and the identity map such that 
$\left(\phi^g(f)\right)_{|\partial C} = f$, and satisfying (3.2). 

Note moreover that $\phi^g(f)$ is $C^{k,\alpha}$-close to the identity if 
$f$ is close to it on $\partial C$ and $g$ is $C^{k-1,\alpha}$-close to 
the Euclidean metric on $C$. This implies that the family $\{\Sigma_{\tau}\}$ 
forms a $C^{k,\alpha}$ foliation of $C$. 
{\endproof}

The metric $g$ in the $x^{\alpha} = (\tau, x^{i})$ coordinates, 
i.e.~$\phi^{*}g$, may be written in lapse/shift form, commonly used in general 
relativity, as
\begin{equation}\label{e3.3}
g = u^{2}d\tau^{2} + g_{ij}(dx^{i} + \sigma^{i}d\tau)(dx^{j} +
 \sigma^{j}d\tau),
\end{equation}
where $u$ is the lapse and $\sigma$ is the shift in the $x$-coordinates 
and $g_{ij}$ is the induced metric on the leaves $S_{\tau} = \{\tau = const\}$. 
A simple computation shows that lapse and shift are related to the metric $g = 
g_{\alpha\beta}^ydy^{\alpha}dy^{\beta}$ in the initial $(y^{\alpha})$ coordinates 
by the equations
\begin{equation}\label{e3.4}
u^{2} + |\sigma|^{2} =
g_{\alpha\beta}^{y}(\partial_{\tau}\phi^{\alpha})(\partial_{\tau}\phi^{\beta}),
\end{equation}
\begin{equation}\label{e3.5}
g_{ij}\sigma^{j} = 
g_{\alpha\beta}^{y}(\partial_{\tau}\phi^{\alpha})(\partial_{i}\phi^{\beta}),
\end{equation}
\begin{equation}\label{e3.6}
g_{ij} = g_{\alpha\beta}^y(\partial_{i}\phi^{\alpha})(\partial_{j}\phi^{\beta}).
\end{equation}
A computation using \eqref{e3.5} shows that 
$|\sigma|^{2} = g^{ij}g_{\alpha\beta}^{y}g_{\mu\nu}^{y}
\partial_{\tau}\phi^{\alpha}\partial_{\tau}
\phi^{\mu}\partial_{i}\phi^{\beta}\partial_{j}\phi^{\nu}$. 
>From $g_{0j} = g_{ij}\sigma^i$ and $g_{00} = u^2 + |\sigma|^2$, 
one may compute $g^{\alpha\beta}$ and, expanding, this yields $g^{00} = u^{-2}$ 
and $\sigma^i = - u^2 g^{0i}$. The unit normal $N$ to the foliation $\Sigma_{\tau}$ 
is given by 
\begin{equation}\label{e3.7}
N = u^{-1}(\partial_{\tau} - \sigma),
\end{equation}
so that, for instance, $g(N,\cdot)=ud\tau$ (this will be useful later on).

\medskip

It is now important to notice that the construction of $H$-harmonic coordinates 
in Lemma 3.2 can be done for any choice of boundary diffeomorphism $f$. We 
shall show that there is a (unique) choice of $f$ close to the identity with 
$f = id$ on $\partial_{0}C = \{0\}\times S^{n-1}$, such that $u$ is 
identically $1$ and the shift $\sigma$ vanishes on the vertical boundary 
$\partial C = I \times S^{n-1}$. 

\medskip

\begin{lemma}\label{l3.3}
For any $k \geq 3$, there exists a $C^{k,\alpha}$ diffeomorphism 
$f: \partial C \rightarrow \partial C$ such that the lapse $u$ and shift 
$\sigma$ of $g$ in \eqref{e3.3} satisfies 
\begin{equation}\label{e3.8}
u = 1, \ {\rm and} \ \sigma = 0,  \ \ {\rm on} \ \ \partial C.
\end{equation}
\end{lemma}

{\bf Proof:}  Consider the operator 
\begin{equation}\label{7}
\Xi: \Met^{k-1,\alpha}(C)\times \Diff_{0}^{k,\alpha}(\partial C) \rightarrow 
C^{k-1,\alpha}(\partial C)\times \prod_{1}^{n}C^{k-1,\alpha}(\partial C),
\end{equation}
$$\Xi(g, f) = (g_{\alpha\beta}^{y}(\partial_{\tau}\phi^{\alpha})
(\partial_{\tau}\phi^{\beta}) - |\sigma|^{2}(\phi), \sigma^{i}(\phi)) ,$$
where $\phi = \phi^g(f)$ is defined above in the proof of Lemma 3.2; recall that 
$\phi|_{\partial C} = f$. More precisely, $\Xi$ is defined in the neighborhoods 
${\mathcal U}$ and ${\mathcal V}$ defined in Lemma 3.2 above. From \eqref{e3.5}, 
one has $\sigma^{i} = g^{ij}g_{\alpha\beta}^{y}\partial_{\tau}
\phi^{\alpha}\partial_{j}\phi^{\beta}$. 
Note that for the map $f = id$ on $\partial C$, and at the metric $g_0 = \delta$, 
one has $\phi^{g_{0}}(id) = id$ and $|\xi|^{2}(id) = 0$, so that 
$\Xi(g_{0}, id) = (1, 0, 0)$. Thus, 
\begin{equation}\label{8}
\Xi(g, id) = (1 + O(\varepsilon), O(\varepsilon))
\end{equation}
where, as already discussed, $\varepsilon$ is positive and may be taken 
as small as needed. We would like to apply the implicit function theorem 
to assert that for any $g \in {\mathcal U}$, where ${\mathcal U}$ is 
sufficiently small, there exists $f = f(g) \in {\mathcal V}$, such that 
\begin{equation}\label{9}
\Xi(g, f(g)) = (1, 0).
\end{equation}
If such $f$ exists, then, for any $g \in {\mathcal U}$, the pair $(g, f)$ defines 
a $C^{k,\alpha}$ diffeomorphism $\phi: C \rightarrow C$ and the resulting 
metric $\phi^{*}g$ satisfies \eqref{e3.8}. Thus it suffices to solve \eqref{9}. 

There is however a loss of one derivative in the map $\Xi$ and its 
derivative in the second variable, as is obvious by looking at its value 
at the metric $g_0=\delta$: 
\begin{equation}\label{10}
(D_2\Xi)_{(g_0,id)}(h) = 
(2\partial_{\tau}h^{0},\partial_{\tau}h^{i} + \partial_{i}h^{0}).
\end{equation}
Thus, we need to use the Nash-Moser inverse function theorem. We use this 
in the form given in [34,\S 6.3], and in particular [34, Thm.6.3.3, Cor.~1, 
Cor.~2]. Following Zehnder's notation, (with $s$ in place of $\sigma$), let 
$X_{s} = \Met^{k-1,\alpha}(C)$, $Y_{s} = \Diff^{k,\alpha}_{0}(\partial C)$, 
and $Z_{s} = C^{k-1,\alpha}(\partial C)\times \prod_{1}^{n}C^{k-1,\alpha}
(\partial C)$, so that $s$ is a linear function of $k+\alpha$. Thus we write 
$X_{s} = \Met^{s + 1 + \varepsilon}(C)$, for some arbitrary but fixed 
$\varepsilon > 0$ (recall we start at $k \geq 2$), $Y_{s} = \Diff^{s + 2 + 
\varepsilon}_{0}(\partial C)$ and $Z_{s} = \prod_{1}^{n}C^{s + 1 + 
\varepsilon}(\partial C)$. We check the hypotheses of Zehnder's theorem:

\smallskip

(H1) When $s = 0$, $\Xi$ is $C^{2}$ in $f$, with uniform bounds in $Y_{0}$. 
This is clearly true. 

(H2) $\Xi$ is Lipschitz in $X_{0}$, also true. 

(H3) $\Xi$ is of order $s = \infty$, with growth $\delta = 1$. This follows 
from
$$||\Xi(g, f)||_{C^{k-2,\alpha}} \leq C(k)(||g||_{C^{k-1,\alpha}} 
+ ||f||_{C^{k,\alpha}}).$$ 

(H4) Existence of right inverse of loss $\gamma = 1$. Let 
$(D_{2}\Xi)_{(g,f)}$ be the derivative of $\Xi$ with respect to the 2nd variable 
$f$ at $(g,f)$. Then varying $f$ in the direction $v$, $f_{s} = f + sv$, it 
is easy to see that the operator $(D_{2}\Xi)_{(g,f)}$ is a 1st order linear 
PDE in $v$, with all coefficients in $C^{k-1,\alpha}$. As in (3.12), the 
boundary $S^{n-1} = \{\tau = 0\}$ is non-characteristic. Hence, for any 
$h \in C^{k-1,\alpha}$, there exists a unique $C^{k-1,\alpha}$ smooth 
solution $v$ to $$D_{2}\Xi_{(g,f)}(v) = h$$ with initial value $v_{0} = 0$ 
on $S^{n-1}$. This gives the existence of an inverse operator $L_{(g,f)}$ 
to $D_{2}\Xi_{(g,f)}$, with a loss of 1-derivative. One has 
$L: Z_{s} \rightarrow Y_{s- 1}$ with $D_{2}\Xi_{(g,f)} \circ L = id$. The 
remaining conditions of (H4) are easily checked to hold. 
It follows then from [34, Cor.~2,p.241] that for any $g$ close 
to $g_{0}$ in $X_{2 + \varepsilon}$ there exists $f \in Y_{1}$, (depending 
continuously on $g$), which satisfies \eqref{9}, (and similarly for higher $s$). 

   This shows that, for any $g \in \Met^{k,\alpha}(C)$ close to $g_{0}$, 
with $k \geq 3$, there exists $f \in C^{k,\alpha'}(\partial C)$, which 
solves \eqref{9}. Pulling back as above gives, for any initial 
$g \in C^{k,\alpha}$, a $C^{k-1,\alpha'}$ metric $\phi^{*}g$ in 
$H$-harmonic coordinates and satisfying \eqref{e3.8}. 
{\endproof}

  For the remainder of the proof, we work in the fixed $H$-harmonic coordinate 
system satisfying \eqref{e3.8}. Next, we derive the form of the Einstein 
equations for the metric $g$ in (\ref{e3.3}). First, the $2^{\rm nd}$ 
fundamental form $A = \frac{1}{2}{\mathcal L}_{N}\, g_{S_{\tau}}$ of 
the leaves $S_{\tau}$ has the form
\begin{equation}\label{e3.17}
A = {\tfrac{1}{2}}u^{-1}(\mathcal{L}_{\partial_{\tau}}g_{S} - 
{\mathcal L}_{\sigma}g_{S}),
\end{equation}
where we have denoted by $g_{S}$ the restriction of $g$ on $S_{\tau}$. More
precisely, and since we shall compute on the $(n+1)$-dimensional manifold 
with tensors living on the $n$-dimensional slices $S_{\tau}$, 
$$ g_{S} = g(\Pi_{S}\cdot,\Pi_{S}\cdot)$$
where $\Pi_{S}$ is the orthogonal projection operator on $S_{\tau}$. Thus,
$g_{S} = g_{ij}(dx^i + \sigma^i d\tau)(dx^j + \sigma^j d\tau)$, as 
in \eqref{e3.3}. Clearly \eqref{e3.17} is the same as
\begin{equation} \label{e3.18}
\mathcal{L}_{\partial_{\tau}}g_{S} = 2uA + {\mathcal L}_{\sigma}g_{S} . 
\end{equation}
A straightforward computation from commuting derivatives gives the Riccati
equation 
\begin{equation} \label{e3.19}
({\mathcal L}_{N}A) = A^{2} - u^{-1}(D^{2}u) - R_{N}, 
\end{equation}
where $R_{N} = g_{S}(R(\cdot,N)N,\cdot)_{|TS\otimes TS}$ and 
$A^{2}$ is the bilinear form associated through $g_{S}$ to the square
of the shape operator of $S_{\tau}$. (The equation (\ref{e3.19}) may also be
derived from the $2^{\rm nd}$ variation formula). Using the fact
that $A$ is tangential, (i.e.~$A(N, \cdot ) = 0$), this gives
\begin{equation} \label{e3.20}
\partial_{\tau}A = -{\mathcal L}_{\sigma}A - D^{2}u + uA^{2} - uR_{N}. 
\end{equation}
Another straightforward calculation via the Gauss equations shows that 
$R_{N} = \Ric_{g} - \Ric_{S_{\tau}} + HA - A^{2}$, which, via 
\eqref{e3.18} and \eqref{e3.20} gives the system of 'evolution' equations for 
$g_{ij}$ and $A = A_{ij}$ on $S_{\tau}$:
\begin{equation} \label{e3.21}
\partial_{\tau}g = 2uA + {\mathcal L}_{\sigma}g_{S}, 
\end{equation}
\begin{equation} \label{e3.22}
\partial_{\tau}A = {\mathcal L}_{\sigma}A - D_{S}^{2}u + u\left( \Ric_{S} 
- \Ric_{g} + 2A^{2} - HA \right). 
\end{equation}
(Up to sign differences, these are the well-known Einstein evolution 
equations in general relativity, cf. \cite{8,32}). Substituting the 
expression of $A$ given by \eqref{e3.21} in (\ref{e3.22}) gives the 
$2^{\rm nd}$-order evolution equation for $g$:
\begin{equation}\label{e3.23}
(\mathcal{L}_{\partial_{\tau}}\mathcal{L}_{\partial_{\tau}} 
+ \mathcal{L}_{\sigma}\mathcal{L}_{\sigma} - 
2\mathcal{L}_{\partial_{\tau}}\mathcal{L}_{\sigma})g_{S} = 
udu(N) A - 2u D_{S}^{2}u + 2u^2\left( \Ric_{S} - \Ric_{g} + 2A^{2} 
- HA \right). 
\end{equation}
We now shift from these intrinsic equations to their expressions in coordinates.
Any tangential $1$-form on $S_{\tau}$ necessarily is of the form
$$ \alpha = \alpha_i(dx^i + \sigma^i d\tau),$$ thus it is enough to work with
the $(i,j)$ components only. 
Using (\ref{e2.7}), (along the slices $S_{\tau}$), one obtains
\begin{equation} \label{e3.24}
 u^{2}\Delta_{S} g_{ij} + \left( (\mathcal{L}_{\partial_{\tau}}\mathcal{L}_{\partial_{\tau}} 
+ \mathcal{L}_{\sigma}\mathcal{L}_{\sigma} 
- 2\mathcal{L}_{\partial_{\tau}}\mathcal{L}_{\sigma})g_{S}\right)_{ij}
= -2u^{2}(\Ric_{g})_{ij}-2u(D_{S}^{2}u)_{ij} + Q_{ij}(g,\partial g), 
\end{equation}
where $Q_{ij}$ is a term involving at most the first order derivatives 
of $g_{\alpha\beta}$ in all $x^{\alpha}$ directions. Now, 
$$ (\mathcal{L}_{\partial_{\tau}}g_{S})_{ij} = \partial_{\tau}g_{ij},\quad 
(\mathcal{L}_{\sigma}g_{S})_{ij} = \sigma^k\partial_k g_{ij} 
+ g_{kj}\partial_i\sigma^k + g_{ik}\partial_j\sigma^k , $$
so that, for Einstein metrics,
\begin{equation} \label{e3.25}
(\partial_{\tau}^{2}+ u^{2}\Delta  
- 2\sigma^k\partial_k\partial_{\tau} + \sigma^k\sigma^l\partial^2_{kl})g_{ij} = 
- 2u(D^{2}u)_{ij} + S_{ij}(g,\partial g) + Q_{ij}(g,\partial g), 
\end{equation}
where $Q_{ij}$ has the same general form as before and $S_{ij}$ contains tangential
first and second derivatives of $\sigma$. 

   The $0i$ and $00$ components of the Ricci curvature in the bulk are given by 
the `constraint' equations along each leaf $S_{\tau}$:
\begin{equation} \label{e3.26} \begin{split}
\delta (A - Hg) & = - \Ric_{g}(N, \cdot  ) = 0, \\
|A|^{2} - H^{2} + R_{S_{\tau}} & = R_{g} - 2\Ric_{g}(N,N) = (n-1)\lambda . 
\end{split}
\end{equation}

\medskip

Next, we derive the equations for the lapse $u$ and shift $\sigma$ 
along the leaves $S_{\tau}$. 

\begin{lemma} \label{l 3.4.}
The lapse $u$ and shift $\sigma$ satisfy the following equations:
\begin{equation} \label{e3.27}
\Delta u + |A|^{2}u + \lambda u = -uN(H) = 
-(\partial_{\tau} - \sigma)H. 
\end{equation}
\begin{equation} \label{e3.28}
\Delta\sigma^{i} = -2u\langle D^{2}x^{i}, A \rangle - u\langle dx^{i}, 
dH \rangle - 2\langle dx^{i}, A(\nabla u) - {\tfrac{1}{2}}Hdu \rangle .
\end{equation}
\end{lemma}

{\bf Proof:} The lapse equation is derived by taking the trace of 
(\ref{e3.19}), and noting that 
$$ tr {\mathcal L}_{N}A = N(H) + 2|A|^{2}.$$ 
For the shift equation, since the functions $x^{i}$ are harmonic on 
$S_{\tau}$, one has
$$\Delta ((x^{i})') + (\Delta')(x^{i}) = 0, $$
where $'$ denotes the Lie derivative with respect to $uN$
and the Laplacian is taken with respect to the induced metric on 
the slices $S_{\tau}$. Moreover $(x^{i})' = - \sigma^{i}$ (see above),
and from standard formulas, cf. \cite[Ch.~1K]{11} for example, one has
$$(\Delta')(x^{i}) = -2\langle D^{2}x^{i}, \delta^{*}uN\rangle 
+ 2\langle dx^{i}, \beta (\delta^{*}uN)\rangle, $$
where all the terms on the right are along $S_{\tau}$ and $\beta$ 
is the Bianchi operator, $\beta(k) = \delta k + \frac{1}{2}dtr k$. 
Thus, $\delta^{*}(uN) = uA$, and the shift components $\sigma^{i}$
 satisfy
$$\Delta\sigma^{i} = -2u\langle D^{2}x^{i}, A\rangle + 2u\langle dx^{i},
\delta A - {\tfrac{1}{2}}dH\rangle + 2\langle dx^{i}, A(\nabla u) - 
{\tfrac{1}{2}}Hdu \rangle.$$
The relation (\ref{e3.28}) then follows from the constraint equation 
(\ref{e3.26}).{\endproof}

Summarizing the work above, the Einstein equations in local $H$-harmonic 
coordinates imply the following system on the data $(g_{ij}, u, \sigma)$:
\begin{equation} \label{e3.31}
(\partial_{\tau}^{2}+ u^{2}\Delta  
- 2\sigma^k\partial_k\partial_{\tau} + \sigma^k\sigma^l\partial^2_{kl})g_{ij} = 
- 2u(D^{2}u)_{ij} + S_{ij}(g_{\alpha\beta}, \partial g_{\alpha\beta}) 
+ Q_{ij}(g_{\alpha\beta}, \partial g_{\alpha\beta}), 
\end{equation}
\begin{equation} \label{e3.32}
\Delta u + |A|^{2}u + \lambda u = dH_0(\sigma). 
\end{equation}
\begin{equation} \label{e3.33}
\Delta\sigma^{i} = -2u\langle D^{2}x^{i}, A \rangle  -
 u\partial_{i}H_0 - 
2\left( A^i_j\nabla^j u - {\tfrac{1}{2}}H_0 \nabla^i u  \right)  , 
\end{equation}
where $H_0$ denotes the mean curvature of the $\{\tau=0\}$-slice $U$.

\begin{remark} \label{r 3.6.}
{\rm  The system (\ref{e3.31})-(\ref{e3.33}) is essentially an elliptic system 
in $(g_{ij}, u, \sigma)$, given that $H = H_{0}$ is prescribed. Thus, assuming 
$u \sim 1$ and $\sigma \sim 0$, the operator $P = \partial_{\tau}^{2}+ 
u^{2}\Delta  - 2\sigma^i\partial_{i}\partial_{\tau} + 
\sigma^k\sigma^l\partial_{kl}^{2}$ is elliptic on $C$ and acts diagonally on
$\{g_{ij}\}$, as is the Laplace operator on the slices $S_{\tau}$ 
acting on $(u, \sigma)$. The system  (\ref{e3.31})-(\ref{e3.33}) is of course 
coupled, but the couplings are all of lower order, i.e.~$1^{\rm st}$ order,
 except for the term $D^{2}u$ in (\ref{e3.31}). However, this term can be 
controlled or estimated by elliptic regularity applied to the lapse equation 
(\ref{e3.32}) (as discussed further below). 
Given the above, it is not difficult to deduce that local $H$-harmonic 
coordinates have the optimal regularity property, i.e.~if $g$ is in 
$C^{m,\alpha}(C)$ in some local coordinate system, then $g$ is in 
$C^{m,\alpha}(C)$ in $H$-harmonic coordinates. Since this will not 
actually be used here, we omit further details of the proof. }
\end{remark}

  Next we show that the lapse and shift, and their $\tau$-derivatives, 
are determined by the tangential metric $g_{S}$ and its $\tau$-derivative. 
\begin{lemma} \label{l 3.5.}
Suppose the metric $g$ is close to the Euclidean metric in the $C^{2,\alpha}$ 
topology. Then in local $H$-harmonic coordinates $(\tau,x^i)$ as defined above, 
the lapse-shift components $(u,\sigma^{i})$ and their derivatives 
$(\partial_{\tau}u,\partial_{\tau}\sigma^{i})$, are uniquely determined 
either by the tangential metric $g_{ij}$ and $2^{\rm nd}$ fundamental form 
$A_{ij}$ on each $S_{\tau}$, or by the tangential metric $g_{ij}$ and its 
time derivatives $\partial_{\tau}g_{ij}$ on each $S_{\tau}$. 
\end{lemma}

{\bf Proof:} The system \eqref{e3.32}-(\ref{e3.33}) is a coupled elliptic 
system in the pair $(u, \sigma)$ on $S_{\tau}$, with boundary values on 
$\partial S_{\tau}$ given by 
\begin{equation} \label{e3.29}
u|_{\partial S_{\tau}} = 1 ,  \ \ \sigma|_{\partial S_{\tau}} = 0. 
\end{equation}
In the $x^{i}$ coordinates, all the coefficients of (\ref{e3.32})-(\ref{e3.33}) 
are bounded in $C^{\alpha}$. Since the metric $g_{ij}$ is close to the flat 
metric in the $C^{2,\alpha}$ topology, it is standard that there is then 
a unique solution to the elliptic boundary value problem 
(\ref{e3.32})-(\ref{e3.33})-(\ref{e3.29}), cf.~[19]. The solution $(u, \sigma)$ 
is uniquely determined by the coefficients $(g_{ij}, A_{ij})$ and the 
terms or coefficients containing derivatives of $H$ and the $x^i$. But 
these are also determined by $(g_{ij}, A_{ij})$. Combining the facts above, 
it follows that $(u, \sigma)$ is uniquely determined by $(g_{ij}, A_{ij})$. 

The second claim is obtained in the same manner: rewrite the equations by
replacing all the occurrences of $A_{ij}$ by its expression in \eqref{e3.17}. 
The equations are then non-linear equations in $(u,\sigma^i)$. Considering 
them as a non-linear operator from $C^{2,\alpha}$ to $C^{\alpha}$ depending 
also on the metric, a simple computation shows that the operator linearized 
at the Euclidean metric is invertible. Invertibility of the non-linear operator 
then follows from the implicit function theorem.

  Next we claim that $\partial_{\tau}g_{0\alpha}$ is also determined by 
$(g_{ij}, A_{ij})$ along $S_{\tau}$. To see this, first note that 
$$ \mathcal{L}_Ng = \mathcal{L}_Ng_{S} + \mathcal{L}_N(g(N,\cdot))\otimes 
g(N,\cdot) + g(N,\cdot)\otimes \mathcal{L}_N(g(N,\cdot)) ,$$
where $\mathcal{L}_Ng_{S}=2A$, $g(N,\cdot) = u d\tau$  and 
$$ \mathcal{L}_N(g(N,\cdot))(\partial_i) = - g(N,[u^{-1}(\partial_{\tau} - 
\sigma,\partial_i])= u^{-1}du(\partial_i), \quad \mathcal{L}_N(g(N,\cdot)))(N) 
= 0.$$
This shows that all components of $\mathcal{L}_Ng$ are determined by $(g_{ij}, 
A_{ij})$, (since $u$ and $\sigma$ are already so determined). Now write 
$\mathcal{L}_Ng(\partial_{\tau}, \partial_{\alpha}) = N(g_{0\alpha}) + 
l_{0\alpha}$. One has $N(g_{0\alpha}) = u^{-1}\partial_{\tau}g_{0\alpha} - 
u^{-1}\partial_{\sigma}g_{0\alpha}$, and the second term is again determined 
by $(g_{ij}, A_{ij})$. Calculating the term $l_{0i}$ above explicitly, 
one easily finds that it also depends only on $(g_{ij}, A_{ij})$, so that
$$\partial_{\tau}g_{0i} = \phi_i  ,$$
is determined by an explicit formula in $g_{ij}$, $A_{ij}$, $u$, $\sigma$ and 
their tangential derivatives, and so implicitly by $g_{ij}$, $A_{ij}$. 
Working now in the same way shows that the same is true for 
$\partial_{\tau}g_{00}$. This completes the proof. 
{\endproof}

\subsection*{Proof of Theorems 3.1 and 1.1.}
Suppose that $g$ and $\widetilde{g}$ are two Einstein metrics on $C$ 
with identical $(\gamma, A)$ on $\partial_{0}C = U$. One may construct 
$H$-harmonic coordinates for each, and via a diffeomorphism identifying these 
coordinates, assume that the resulting pair of metrics $g$ and $\widetilde g$ 
have fixed $H$-harmonic coordinates $(\tau, x^{i})$, and both metrics satisfy 
the system (\ref{e3.31})-(\ref{e3.33}). Let
\begin{equation} \label{e3.34}
h = h_{ij} = \widetilde g_{ij}- g_{ij}. 
\end{equation}
One then takes the difference of both equations (\ref{e3.31}) and freezes the 
coefficients at $g$ to obtain a linear equation in $h$. Thus, for example, 
$\Delta_{\widetilde g}\widetilde g_{ij} - \Delta_{g}g_{ij} = 
\Delta_{g}(h_{ij}) - (g^{ab} - \widetilde g^{ab})
\partial_{a}\partial_{b}\widetilde g_{ij}$. The second 
term here is of zero order, (rational), in the difference $h$, with 
coefficients depending on two derivatives of $\widetilde g$. Carrying out 
the same procedure on the remaining terms in \eqref{e3.31} gives the equation
$$(\partial_{\tau}^{2}+ u^{2}\Delta  - 
2\partial_{\sigma}\partial_{\tau} + \partial_{\sigma}^{2})h_{ij} = - 
2(\widetilde u(\widetilde D^{2}\widetilde u)_{ij} - u(D^{2}u)_{ij}) + 
\mathcal{Q}_{ij}(h_{\alpha\beta},\partial_{\mu} h_{\alpha\beta}),$$
where we have denoted $\partial_{\sigma}=\sigma^k\partial_k$, 
$\partial_{\sigma}^2=\sigma^k\sigma^l\partial_{kl}$, and $\mathcal{Q}$ is
a term depending on two derivatives of the background $(g,\widetilde g)$
and linear in its arguments, whose precise value may change from line to line.
Similarly, 
$\widetilde D^{2}\widetilde u - D^{2}u = D^{2}v + (\widetilde D^{2} -
D^{2})\widetilde u$, where $v = \widetilde u - u$ and the second term is 
of the form $Q$ above. Hence, 
\begin{equation} \label{e3.35}
(\partial_{\tau}^{2}+ u^{2}\Delta  - 
2\partial_{\sigma}\partial_{\tau} + \partial_{\sigma}^{2})h_{ij} = - 
2u(D^{2}v)_{ij} + \mathcal{Q}_{ij}(h_{\alpha\beta},\partial_{\mu}
 h_{\alpha\beta}), 
\end{equation}
Note that since we have linearized, $\mathcal{Q}$ depends linearly on
$h_{\alpha\beta}$ and $\partial_{\mu} h_{\alpha\beta}$, with nonlinear 
coefficients depending on $\widetilde g$ and $g$. 

Next we use the lapse and shift equations (\ref{e3.32})-(\ref{e3.33}) 
to estimate the differences $v = \widetilde u - u$ and 
$\chi = \widetilde \sigma - \sigma$. Thus, as before, 
$\Delta_{\widetilde g}\widetilde u 
- \Delta_{g}u = \Delta_{g}v + D^{2}_{h}(\widetilde u)$, 
where $D^{2}_{h}$ is a $2^{\rm nd}$-order differential operator on $\widetilde u$ 
with coefficients depending on the difference $h$, to first order. The remaining 
terms in (\ref{e3.32})-(\ref{e3.33}) can all be treated in the same way, using 
\eqref{e3.17} to replace occurences of $A_{ij}$ by derivatives in $\tau$ and 
$\sigma$. Taking the difference, it then follows from (\ref{e3.32}) and 
(\ref{e3.33}) that 
\begin{equation} \label{e3.36}\begin{cases}
\Delta v + |A|^2 v +\lambda v  = Q(h_{ij}, \partial_{\mu}h_{ij}) \\
\Delta \chi^i +2v\langle D^2x^i,A\rangle + v \partial_iH_0 
+2 (A_j^i-\frac{1}{2}H\delta_j^i)g^{jk}\partial_ku
 = Q^{i}(h_{kl}, \partial_{\mu}h_{kl})  
\end{cases} 
\end{equation}
where the terms $Q$ are linear in the arguments and their 
coefficients depend on one derivatives of 
$\widetilde g$. Note also that the zeroth- and first-order terms in $(v,\chi)$ 
aresmall if the metric is close to the Euclidean metric. Thus, the left-hand 
side operators are invertible with $v = 0$ and $\chi=0$ on $\partial S_{\tau}$, 
and elliptic regularity applied to the system \eqref{e3.36} then gives
\begin{equation} \label{e3.37}
||v||_{L_{x}^{2,2}} \leq  C(\widetilde g, g) ||h_{ij}||_{L_{(\tau ,x)}^{1,2}} , 
\end{equation}
and
\begin{equation} \label{e3.38}
||\chi||_{L_{x}^{2,2}} \leq  C(\widetilde g, g) 
||h_{ij}||_{L_{(\tau ,x)}^{1,2}} . 
\end{equation}
It follows from (\ref{e3.35}) and (\ref{e3.37})-(\ref{e3.38}) that
\begin{equation} \label{e3.39}
||P(h_{ij})||_{L_{x}^{2}} \leq  C(\widetilde{g}_{\alpha\beta}, g_{\alpha\beta})
\, ||h_{\alpha\beta}||_{L_{(\tau ,x)}^{1,2}}, 
\end{equation}
where $P$ is given in Remark \ref{r 3.6.}. 

   Now by applying Lemma \ref{l 3.5.} to $(u, \sigma)$ and $(\widetilde u, 
\widetilde \sigma)$ and taking the difference as above, it follows that $v$ 
and $\chi$, as well as $\partial_{\tau}v$ and $\partial_{\tau}\chi$ are 
given by a linear expression in $h_{ij}$ and its first derivatives 
(in every direction). Hence, \eqref{e3.39} becomes
\begin{equation}\label{e3.40} 
||P(h_{ij})||_{L_{x}^{2}} \leq  C(\widetilde g, g)
||h_{ij}||_{L_{(\tau ,x)}^{1,2}}. 
\end{equation}

   We are now in position to apply the Calder\'on unique continuation 
theorem [14]. Thus, the operator $P$ is elliptic and diagonal, and the 
Cauchy data for $P$ vanish at $U$, i.e.
\begin{equation} \label{e3.41}
h = \partial_{\tau}h = 0 \ \ {\rm at} \ \  U.
\end{equation}
We claim that $P$ satisfies the hypotheses of the Calderon unique 
continuation theorem [14]. Following [14], decompose the symbol of $P$ as 
\begin{equation} \label{e3.42}
A_{2}(\tau ,x,\xi) = (u^{2}g^{kl}\xi_{k}\xi_{l} -
2\sigma^{k}\sigma^{l}\xi_{k}\xi_{l})I, 
\end{equation}
$$A_{1}(\tau ,x,\xi) = \sigma^{k}\xi_{k}I, $$
where $I$ is the $N\times N$ identity matrix, $N = \frac{1}{2}n(n+1)$, 
equal to the cardinality of $\{ij\}$. Setting $|\xi|^{2} = 1$, \eqref{e3.42} 
becomes
$$A_{2}(\tau ,x,\xi) = (u^{2} - 
2\sigma^{k}\sigma^{l}\xi_{k}\xi_{l})I,$$
$$A_{1}(\tau ,x,\xi) = \sigma^{k}\xi_{k}I.$$
Now form the matrix
\begin{equation}\label{e3.43}
M = \left(
\begin{array}{cc}
0  &  -I \\
A_{2} & A_{1}
\end{array}
\right)
\end{equation}
The matrices $A_{1}$ and $A_{2}$ are diagonal, and it is then easy to 
see that $M$ is diagonalizable, i.e.~has a basis of eigenvectors over 
${\mathbb C}$. This implies that $M$ satisfies the hypotheses of 
[14, Thm.~11(iii)], cf.~also [14, Thm.~4]. The bound \eqref{e3.40} is 
substituted in the basic Carleman estimate of [14, Thm.~6], cf.~also 
[29, (6.1)], showing that $h_{ij}$ satisfies the unique continuation 
property. It follows from \eqref{e3.41} and the Calder\'on unique 
continuation theorem that 
$$h_{ij} = \widetilde g_{ij} - g_{ij} = 0, $$
in an open neighborhood $\Omega \subset C$. 

By Lemma \ref{l 3.5.} once again, this implies $h_{\alpha\beta} = 0$, i.e.
$$\widetilde g_{\alpha\beta} = g_{\alpha\beta}, $$
in $\Omega$, so that $\widetilde{g}$ is isometric to $g$ in $\Omega$. 
By construction, the isometry from $\widetilde{g}$ to $g$ equals the identity 
on $U$. This shows that the metric $g$ is uniquely determined in $\Omega$,
up to isometry, by the abstract Cauchy data on $U$. Since Einstein 
metrics are real-analytic in the interior in harmonic coordinates, a
standard analytic continuation argument, (cf.~[25] for instance), then 
implies that $g$ is unique up to isometry everywhere in $C$. This 
completes the proof of Theorem 3.1. 

   In the context of Theorem 1.1, the same analytic continuation argument 
shows that a pair of Einstein metrics $(M_{i}, g_{i})$, $i = 0,1$, whose 
Cauchy data agree on a common open set $U$ of $\partial M_{i}$ are 
everywhere locally isometric, i.e.~they become isometric in suitable 
covering spaces, modulo restriction or extension of the domain, as 
discussed following Theorem 1.1. This then also completes the proof 
of Theorem 1.1.
{\endproof}

As an illustration, suppose $(M_{1}, g_{1})$ and $(M_{2}, g_{2})$ are a 
pair of Einstein metrics on compact manifolds-with-boundary and the Cauchy 
data for $g_{1}$ and $g_{2}$ agree on an open set $U$ of the boundary. 
Suppose $M_{i}$ are connected and the topological condition (1.5) holds 
for each $M_{i}$. Then, modulo isometry, either $M_{1} \subset M_{2}$, 
$M_{2} \subset M_{1}$, or $M_{i}$ are subdomains in a larger Einstein 
manifold $M_{3} = M_{1}\cup M_{2}$. 

\medskip

 We conclude this section with a discussion of generalizations of 
Theorem 1.1. First, one might consider the unique continuation problem 
for 
\begin{equation} \label{e3.45}
\Ric_{g} = T, 
\end{equation}
where $T$ is a fixed symmetric bilinear form on $M$, at least
 $C^{\alpha}$ up to $\bar M$. However, this problem is not natural, 
in that is not covariant under changes by diffeomorphism. For metrics 
alone, the Einstein equation (1.2) is the only equation covariant under 
diffeomorphisms which involves at most the $2^{\rm nd}$ derivatives of 
the metric.  Nevertheless, the proof of Theorem 1.1 shows that if 
$\widetilde g$ and $g$ are two solutions of \eqref{e3.45} which have 
common $H$-harmonic coordinates near (a portion of) $\partial M$ on which 
$(\gamma, A) = (\widetilde \gamma, \widetilde A)$, then $\widetilde g$ is 
isometric to $g$ near (a portion of) $\partial M$.

   Instead, it is more natural to consider the Einstein equation 
coupled (covariantly) to other fields $\chi$ besides the metric; such 
equations arise naturally in many areas of physics. For example, $\chi$ 
may be a function on $M$, i.e.~a scalar field, or $\chi$ may be a 
connection 1-form (gauge field) on a bundle over $M$. We assume that 
the field(s) $\chi$ arise via a diffeomorphism-invariant Lagrangian 
${\mathcal L} = {\mathcal L}(g, \chi)$, depending on $\chi$ and its 
first derivatives in local coordinates, and that $\chi$ satisfies 
field equations, i.e.~Euler-Lagrange equations, coupled to the metric. 
For example, for a free massive scalar field, the equation is the 
eigenfunction equation
\begin{equation} \label{e3.46}
\Delta_{g}\chi  = \mu\chi ,
\end{equation}
while for a connection 1-form, the equations are the Yang-Mills 
equations, (or Maxwell equations when the bundle is a $U(1)$ bundle):
\begin{equation} \label{e3.47}
dF = d^{*}F = 0,
\end{equation}
where $F$ is the curvature of the connection $\chi$. Associated to 
such fields is the stress-energy tensor $T = T_{\mu\nu}$; this is a 
symmetric bilinear form obtained by varying the Lagrangian for $\chi$ 
with respect to the metric, cf.~[22] for example. For the free massive 
scalar field $\chi$ above, one has
$$T = d\chi\cdot  d\chi  - \tfrac{1}{2}(|d\chi|^{2} + \mu\chi^{2})g, $$
while for a connection 1-form 
$$T = F \cdot F - \tfrac{1}{4}|F|^{2}g, $$
where $(F \cdot F)_{ab} = F_{ac}F_{bd}g^{cd}$. 

 When the part of the Lagrangian involving the metric to $2^{\rm nd}$
 order 
only contains the scalar curvature, i.e.~the Einstein-Hilbert action,
 the 
resulting coupled Euler-Lagrange equations for the system $(g, \chi)$ 
are
\begin{equation} \label{e3.48}
\Ric_{g} - \frac{R}{2}g = T, \ \ E_{g}(\chi) = 0. 
\end{equation}
By taking the trace, this can be rewritten as
\begin{equation} \label{e3.49}
\Ric_{g} = \hat T = T - \frac{1}{n-1}tr_{g}T,  \ \ E_{g}(\chi) = 0. 
\end{equation}

 Here we assume $E_{g}(\chi)$ is a $2^{\rm nd}$ order elliptic system
 for 
$\chi$, with coefficients depending on $g$, as in \eqref{e3.46} or 
\eqref{e3.47}, (the latter viewed as an equation for the connection). 
In case the field(s) $\chi$ have an internal symmetry group, as in 
the case of gauge fields, this will require a particular choice of 
gauge for $\chi$ in which the Euler-Lagrange equations become an 
elliptic system in $\chi$. It is also assumed that solutions $\chi$ 
of $E_{g}(\chi) = 0$ satisfy the unique continuation property; for
 instance $E_{g}$ satisfies 
the hypotheses of the Calder\'on theorem [14]. Theorem 1.1 now easily
 extends to cover 
\eqref{e3.48} or \eqref{e3.49}. 

\begin{proposition} \label{p3.7}
Let $M$ be a compact manifold with boundary $\partial M$. Then 
$C^{3,\alpha}$ solutions $(g, \chi)$ of \eqref{e3.48} on $\bar M$ 
are uniquely determined, up to local isometry and inclusion, by 
the Cauchy data $(\gamma, A)$ of $g$ and the Cauchy data 
$(\chi, \partial_{t}\chi)$ on an open set $U \subset \partial M$. 
\end{proposition}

{\bf Proof:}
 The proof is the same as the proof of Theorem 1.1. Briefly, via a
 suitable 
diffeomorphism equal to the identity on $\partial M$, one brings a pair
 of solutions 
of \eqref{e3.48} with common Cauchy data into a fixed system of $H$-harmonic
 coordinates 
for each metric. As before, one then applies Calder\'on uniqueness to
 the resulting 
system \eqref{e3.48} in the difference of the metrics and fields. Further
 details are left 
to the reader. 
{\endproof}

\section{Proof of Theorem 1.2.}
\setcounter{equation}{0}

 Let $g$ be a conformally compact metric on a compact $(n+1)$-manifold
 $M$ with boundary 
which has a $C^{2}$ geodesic 
compactification
\begin{equation} \label{e4.1}
\bar g = t^{2}g, 
\end{equation}
where $t(x) = dist_{\bar g}(x, \partial M)$. By the Gauss Lemma, one
 has the splitting
\begin{equation} \label{e4.2}
\bar g = dt^{2} + g_{t}, 
\end{equation}
near $\partial M,$ where $g_{t}$ is a curve of metrics on $\partial M$ 
with $g_{0} = \gamma$ the boundary metric. The curve $g_{t}$ is 
obtained by taking the induced metric the level sets $S(t)$ of $t$, 
and pulling back by the flow of $N = \bar \nabla t$. Note that if 
$r = -\log t$, then $g = dr^{2} + t^{-2}g_{t}$, so the integral curves 
of $\nabla r$ with respect to $g$ are also geodesics. Each choice of
 boundary 
metric $\gamma \in [\gamma]$ determines a unique geodesic defining 
function $t$.

 Now suppose $g$ is Einstein, so that (1.4) holds and suppose for the
 moment that 
$g$ is $C^{2}$ conformally compact with $C^{\infty}$ smooth boundary
 metric $\gamma$. 
Then the boundary regularity result of [16] implies that $\bar g$ is
 $C^{\infty}$ 
smooth when $n$ is odd, and is $C^{\infty}$ polyhomogeneous when $n$ is
 even. Hence, 
the curve $g_{t}$ has a Taylor-type series in $t$, called the
 Fefferman-Graham 
expansion [18]. The exact form of the expansion depends on whether $n$
 is odd 
or even. If $n$ is odd, one has a power series expansion 
\begin{equation} \label{e4.3}
g_{t} \sim g_{(0)} + t^{2}g_{(2)} + \cdots + t^{n-1}g_{(n-1)} +
 t^{n}g_{(n)} 
+ \cdots , 
\end{equation}
while if $n$ is even, the series is polyhomogeneous,
\begin{equation} \label{e4.4}
g_{t} \sim g_{(0)} + t^{2}g_{(2)} + \cdots + t^{n}g_{(n)} + t^{n}\log t
 \ {\mathcal H} + 
\cdots . 
\end{equation}
In both cases, this expansion is even in powers of $t$, up to $t^{n}$.
 It is important 
to observe that the coefficients $g_{(2k)}$, $k \leq [n/2]$, as well as
 the coefficient 
${\mathcal H}$ when $n$ is even, are explicitly determined by the
 boundary metric $\gamma  = 
g_{(0)}$ and the Einstein condition (1.4), cf.~[18], [20]. For $n$
 even, the series (4.4) 
has terms of the form $t^{n+k}(\log t)^{m}$. 

  For any $n$, the divergence and trace (with respect to $g_{(0)} =
 \gamma$) of $g_{(n)}$ 
are determined by the boundary metric $\gamma$; in fact there is a
 symmetric bilinear form 
$r_{(n)}$ and scalar function $a_{(n)}$, both depending only on
 $\gamma$ and its derivatives 
up to order $n$, such that 
\begin{equation} \label{e4.5}
\delta_{\gamma}(g_{(n)} + r_{(n)}) = 0, \ \ {\rm and} \ \
 tr_{\gamma}(g_{(n)} + r_{(n)}) = 
a_{(n)}.
\end{equation}
For $n$ odd, $r_{(n)} = a_{(n)} = 0$. (The divergence-free tensor
 $g_{(n)} + r_{(n)}$ is 
closely related to the stress-energy of a conformal field theory on
 $(\partial M, \gamma)$, 
cf.~[17]). The relations (4.5) will be discussed further in \S 5. 

  However, beyond the relations (4.5), the term $g_{(n)}$ is not
 determined by $g_{(0)}$; 
it depends on the ``global'' structure of the metric $g$. The higher
 order coefficients 
$g_{(k)}$ of $t^{k}$ and coefficients $h_{(km)}$ of $t^{n+k}(\log
 t)^{m}$, are then 
determined by $g_{(0)}$ and $g_{(n)}$ via the Einstein equations. The
 equations (4.5) 
are constraint equations, and arise from the Gauss-Codazzi and Gauss
 and Riccati equations 
on the level sets $S(t) = \{x: t(x) = t\}$ in the limit $t \rightarrow
 0$; this is also 
discussed further in \S 5. 

   In analogy to the situation in \S 3, the term $g_{(n)}$ corresponds
 to the $2^{\rm nd}$ 
fundamental form $A$ of the boundary, in that, modulo the constraints
 (4.5), it is freely 
specifiable as Cauchy data, and is the only such term depending on
 normal derivatives of 
the boundary metric.

\medskip

Suppose now $g_{0}$ and $g_{1}$ are two solutions of 
\begin{equation} \label{e4.6}
\Ric_{g} + ng = 0, 
\end{equation}
with the same $C^{\infty}$ conformal infinity $[\gamma]$. Then there
 exist 
geodesic defining functions $t_{k}$ such that $\bar g_{k} =
 (t_{k})^{2}g_{k}$ 
have a common boundary metric $\gamma \in [\gamma]$, and both metrics
 are 
defined for $t_{k} \leq  \varepsilon$, for some $\varepsilon > 0$. 

  The hypotheses of Theorem 1.2, together with the discussion above
 concerning 
(4.3) and (4.4), then imply that 
\begin{equation} \label{e4.7}
|g_{1} - g_{0}| = o(e^{-nr}) = o(t^{n}), 
\end{equation}
where the norm is taken with respect to $g_{1}$, (or $g_{0}$). 

\medskip

  Given this background, we prove the following more general version of
 Theorem 1.2, 
analogous to Theorem 3.1. Let $\Omega$ be a domain diffeomorphic to
 $I\times 
B^{n}$, where $B^{n}$ is a ball in ${\mathbb R}^{n}$ with boundary $U =
 \partial 
\Omega$ diffeomorphic to a ball in ${\mathbb R}^{n} \simeq \{0\}\times 
{\mathbb R}^{n}$. 

\begin{theorem} \label{t4.1}
Let $g_{0}$ and $g_{1}$ be a pair of conformally compact Einstein
 metrics on 
a domain $\Omega$ as above. Suppose $g_{0}$ and $g_{1}$ have
 $C^{3,\alpha}$ 
geodesic compactifications, and \eqref{e4.7} holds in $\Omega$. 

  Then $(\Omega, g_{0})$ is isometric to $(\Omega, g_{1})$, by an
 isometry equal 
to the identity on $\partial \Omega$. Hence, if $(M_{0}, g_{0})$ and
 $(M_{1}, g_{1})$ 
are conformally compact Einstein metrics on compact manifolds with
 boundary, and 
\eqref{e4.7} holds on some open domain $\Omega$ in $M_{0}$ and $M_{1}$, 
then the manifolds 
$M_{0}$ and $M_{1}$ are diffeomorphic in some covering space of each
 and the lifted 
metrics $g_{0}$ and $g_{1}$ are isometric.
\end{theorem}

 The proof of Theorem 4.1 is very similar to that of Theorem 3.1. For 
clarity, we first prove the result in case the metrics $g_{i}$, $i = 
0,1$, have a common $C^{\infty}$ boundary metric $\gamma$ and then show 
how the proof can be extended to cover the more general case of metrics 
with less regularity. 

 By applying a diffeomorphism if necessary, one may assume that the 
metrics $g_{i}$ have a common geodesic defining function $t$ defined 
near $\partial\Omega$ and common geodesic boundary coordinates. By
 [16], 
the geodesically compactified metrics $\bar g_{i} = t^{2}g_{i}$ are 
$C^{\infty}$ polyhomogeneous and extend $C^{\infty}$ polyhomogeneously 
to $\partial\Omega$. It follows from the discussion of the Fefferman-Graham 
expansion following (4.5) that $g_{0}$ and $g_{1}$ agree to infinite order 
at $\partial U$, i.e.
\begin{equation} \label{e4.8}
k = g_{1} - g_{0} = O(t^{\nu}), 
\end{equation}
for any $\nu  < \infty$. Of course $k_{0\alpha} = 0$. 

 For the rest of the proof, we work in the setting of the compactified 
metrics $\bar g_{i}$. As in the proof of Theorem 3.1, we assume 
that the domain $\Omega$, now denoted $C$, is sufficiently small so that 
$(C, \bar g_{i})$ is close to the flat metric on the standard cylinder 
$C = I\times B^{n}$, with $\bar A = 0$ on $U = \partial_{0}C$. 
(Note that $g_{(1)} = 0$ in (4.3)-(4.4)). In particular, near 
$\partial_{0}C$, $\bar H = O(t)$. One may construct a foliation 
$S_{\tau}$ with $\bar H_{S_{\tau}} = 0$, together with corresponding 
$H$-harmonic coordinates $(\tau , x^{i})$, exactly as in Lemmas 3.2 
and 3.3, and satisfying the boundary conditions \eqref{e3.8}. All of 
the analysis carried out in \S 3 carries over to this situation with 
only a single difference. Namely, for the term $\Ric_{g}$ in \eqref{e3.23} 
or \eqref{e3.24}, one now no longer has $\Ric_{g} = \lambda g$, but instead 
the Ricci curvature $\bar \Ric$ of the compactified metric $\bar g$. 
Using the facts that $\Ric_{g} = -ng$ and the compactification 
$\bar g$ is geodesic, standard formulas for the behavior of Ricci
 curvature under conformal change give
\begin{equation} \label{e4.9}
\bar \Ric = -(n-1)t^{-1}\bar D^{2}t  - t^{-1}\bar \Delta t \bar g. 
\end{equation}
One has $\bar D^{2}t = {\mathcal L}_{\nabla t}\bar g = O(t)$. 
If $(t, y^{i})$ are geodesic boundary coordinates, then $\partial_{x^{i}} 
= \sum(1-\varepsilon (\tau))\partial_{y^{j}} + \varepsilon (\tau)\nabla t$, 
where $\varepsilon(\tau) = O(\tau)$. Similarly, $\tau /t = 1 + 
\varepsilon(\tau)$. (The specific form of $\varepsilon (\tau)$ 
of course differs in each occurence above, but this is insignificant). 
Since $\bar D^{2}t$ vanishes on $\nabla t$, it follows from (4.9) 
that in the $x^{i}$ coordinates on $S_{\tau}$,
\begin{equation} \label{e4.10}
\bar \Ric_{ij} = -(n-1)(1-\varepsilon)^{2}t^{-1}({\mathcal L}_{\nabla t}
\bar g)_{ij}  - (1-\varepsilon)^{2}t^{-1}(\bar \Delta t)\bar g_{ij} + 
\varepsilon t^{-1}(\bar \Delta t)q_{ij}, 
\end{equation}
where $q_{ij}$ depends only on $\bar g_{0\alpha}$ to zero-order. Next 
$({\mathcal L}_{\nabla t}\bar g) = (1 - \varepsilon)\partial_{\tau}\bar
 g + 
\varepsilon (\tau)\partial_{x^{\alpha}}\bar g$ and similarly for the 
Laplace term in (4.10). Substituting (4.10) in \eqref{e3.24}, it follows 
that the analogue of \eqref{e3.25} in this context is the 'evolution equation'
\begin{equation} \label{e4.11}
\tau^{2}(\partial_{\tau}^{2}+ u^{2}\Delta  - 
2\partial_{\sigma}\partial_{\tau} + \partial_{\sigma}^{2})g_{ij} = - 
2\tau^{2}u(D^{2}u)_{ij} + S_{ij}(g,\tau\partial g) + Q_{ij}(g,\tau\partial g), 
\end{equation}
where $S_{ij}$ and $Q_{ij}$ have the same meaning as before. Here and below, 
we drop the bar from the notation. 

 The lapse $u$ and shift $\sigma$ satisfy essentially the same 
equations as before, namely
\begin{equation} \label{e4.12}
\Delta u + |A|^{2}u - (t^{-1}\Delta t)u = 0, 
\end{equation}
\begin{equation} \label{e4.13}
\Delta\sigma^{i} = -2u\langle D^{2}x^{i}, A \rangle - 2\langle dx^{i},
 A(\nabla u) \rangle . 
\end{equation}
Comparing with (\ref{e3.27})-(\ref{e3.28}), one has here $H = 0$, with the 
$\lambda$ term in replaced by $-t^{-1}\Delta t$. Lemma 3.6 holds 
as before, since $t^{-1}\Delta t$ is smooth up to $\partial_{0}C$. 

 One now proceeds just as in the proof of Theorem 3.1, taking the 
difference of the equation (4.11) to obtain a linear equation on $h = 
\widetilde g - g$; (recall that the bars have been removed from the 
notation). Note that by (4.8), together with elliptic regularity 
applied to (4.12)-(4.13), as in the proof of Lemma 3.6, one has
\begin{equation} \label{e4.14}
h_{\alpha\beta} = O(t^{\nu}), 
\end{equation}
for all $\nu  < \infty$.  The estimates \eqref{e3.37}-\eqref{e3.39} 
and \eqref{e3.40} hold as before.

 Let $P(h_{ij}) = \tau^{2}(\partial_{\tau}^{2}+ u^{2}\Delta  - 
2\partial_{\sigma}\partial_{\tau} + \partial_{\sigma}^{2})$. Then 
$P$ is a fully degenerate $2^{\rm nd}$ order elliptic operator, with 
smooth coefficients, and one has
$$||P(h_{ij})||_{L_{x}^{2}} \leq  C||h_{ij}||_{L_{\tau,x}^{1,2}},$$
where the $1^{\rm st}$ order derivatives on the right are of the form 
$\tau\partial$. Further, by (4.14), $h$ vanishes to infinite order at 
$\partial_{0}C$. It then follows from a unique continuation theorem 
of Mazzeo, [27, Thm.~14], that
$$h_{ij} = 0 $$
in $\Omega \subset C$. The vanishing of $h = h_{\alpha\beta}$ in 
$C$ then follows as before in the proof of Theorem 3.1. 

  Next suppose $g_{0}$ and $g_{1}$ have only a $C^{3,\alpha}$ geodesic 
compactification with a common boundary metric $\gamma$, but that 
(4.7) holds. All of the arguments above remain valid, except the 
infinite order vanishing property (4.8), and the corresponding (4.14), 
which are replaced by the statements $k = o(t^{n})$ and $h = o(t^{n})$ 
respectively. The unique continuation result in [27] per se, requires 
the infinite order decay (4.14). Thus, it suffices to show that (4.14) 
does in fact hold. 

 To do this, we first show that $k = O(t^{\nu})$ weakly, for all $\nu  
< \infty$. This will imply $h = O(t^{\nu})$ weakly, and the strong or 
pointwise decay (4.14) then follows from elliptic regularity.

  In geodesic boundary coordinates, the geodesic compactification of a 
conformally compact Einstein metric satisfies the equation
\begin{equation} \label{e4.15}
t\ddot g - (n-1)\dot g - 2Hg^{T}-2t\Ric_{S(t)} + tH\dot g - t(\dot
 g)^{2} = 0,
\end{equation}
where $\dot g$ is the Lie derivative of $g$ with respect to $\nabla t$, 
cf.~[18] or [21]. Thus $\dot g = 2A$, where $A$ is the $2^{\rm nd}$ 
fundamental form of the level set $S(t)$ of $t$, (with respect to the 
inward normal). Also $H = tr A$, $T$ denotes restriction or projection 
onto $S(t)$ and $\Ric_{S(t)}$ is the intrinsic Ricci curvature of
 $S(t)$. 
(The equation (4.15) may be derived from \eqref{e3.22} by setting $u = 1$ 
and $\sigma = 0$). We recall, as above, that the bar has been removed 
from the notation. 

  As above, the metrics $g_{0}$ and $g_{1}$ are assumed to have a fixed 
geodesic defining function $t$ with common boundary metric $\gamma$ and 
common geodesic boundary coordinates. Taking the difference of the 
equation (4.15) evaluated on $g_{1}$ and $g_{0}$ gives the following 
equation for $k = g_{1} - g_{0}$ as in (4.8):
\begin{equation} \label{e4.16}
t\ddot k - (n-1)\dot k = tr(\dot k) g_{0}^{T} + 2t(\Ric_{S(t)}^{1} -
 \Ric_{S(t)}^{0}) 
+ O(t)k + O(t^{2})\dot k,
\end{equation}
where $O(t^{k})$ denotes terms of order $t^{k}$ with coefficients
 depending 
smoothly on $g_{0}$. One has $\Ric_{S(t)} = D_{x}^{2}(g_{ij})$ is a
 $2^{\rm nd}$ 
order operator on $g_{ij}$, so that (4.16) gives
\begin{equation} \label{e4.17}
t\ddot k - (n-1)\dot k = tr(\dot k) g_{0}^{T} + 2tD_{x}^{2}(k) + O(t)k
 + 
O(t^{2})\dot k,
\end{equation}

   The (positive) indicial root of the trace-free part of (4.16) or
 (4.17) 
is $n$, in that the formal power series solution of (4.17) has
 undetermined 
coefficient at order $t^{n}$, as in the Fefferman-Graham expansion
 (4.3)-(4.4). 
The hypothesis (4.7) implies that 
\begin{equation} \label{e4.18}
k = o(t^{n}), 
\end{equation}
so that this $n^{\rm th}$ order coefficient vanishes. However, taking
 the 
trace of (4.17) gives
$$t \, tr\ddot k - (2n-1)tr \dot k = tr (O(t)k + O(t^{2}\dot k)) + 
2t \, tr (D_{x}^{2}(k)),$$
which has indicial root $2n$. To see that $tr k$ is in fact formally 
determined at order $2n$, one uses the trace of the Riccati equation 
\eqref{e3.19}, (with $u = 1$ and $\sigma = 0$), which gives
\begin{equation}\label{e4.19}
\dot H  + |A|^{2} = -\Ric(T,T).
\end{equation}
Via (4.9), this is easily seen to be equivalent to 
$$t\dot H - H = -t|A|^{2}.$$
This holds for each compactified metric $g_{1}$ and $g_{0}$, and so
 taking 
the difference, and computing as in (4.16)-(4.17) gives the equation
\begin{equation} \label{e4.20}
t\frac{d^{2}}{dt^{2}}(tr k)  - \frac{d}{dt}(tr k) = O(t)k +
 O(t^{2})\dot k .
\end{equation}

  The positive indicial root of (4.20) is 2, and by (4.7), the
 $O(t^{2})$ 
component of the formal expansion of $tr k$ vanishes. Similarly, the
 trace-free 
part $k_{0}$ of $k$ satisfies the equation
\begin{equation} \label{e4.21}
t\ddot k_{0} -(n-1)\dot k_{0} = 2t(D_{x}^{2}(k))_{0} + [O(t)k]_{0} + 
[O(t^{2})\dot k]_{0} , 
\end{equation}
with indicial root $n$. As in [18], by repeated differentiation of
 (4.20) 
and (4.21) it follows from (4.7) that the formal expansion of $k$
 vanishes. 

  Next we show that (4.8) holds weakly. 

\begin{lemma} \label{l4.2}
Suppose $k = o(t^{n})$ weakly, in that, with respect to the
 compactified metric 
$(S(t), g)$, ($g = g_{0}$), 
\begin{equation} \label{e4.22}
\int_{S(t)}\langle k, \phi \rangle = o(t^{n}), \ {\rm as} \ t
 \rightarrow 0 ,
\end{equation}
where $\phi$ is any symmetric bilinear form, $C^{\infty}$ smooth up to 
$U = \partial_{0}C$ and vanishing to infinite order on $\partial C$. 
Then 
\begin{equation} \label{e4.23}
k = o(t^{\nu}), \ {\rm weakly},
\end{equation}
for any $\nu < \infty$, i.e.~{\rm (4.22)} holds, with $\nu$ in place of
 $n$. 
\end{lemma}

{\bf Proof:} 
Here smoothness is measured with respect to the given geodesic
 coordinates 
$(t, x^{i})$ covering $C$. The proof proceeds by induction,
 starting at 
the initial level $n$. As above, the trace-free and pure trace cases
 are treated 
separately, and so we assume in the following first that $\phi$ is
 trace-free. 
Pair $k$ with $\phi$ and integrate (4.17) over the level sets $S(t)$ to
 obtain
\begin{equation} \label{e4.24}
t\int_{S(t)}\langle \ddot k, \phi \rangle - (n-1)\int_{S(t)}\langle
 \dot k, 
\phi \rangle  = t\int_{S(t)}\langle k, P_{2}(\phi) \rangle + 
\int_{S(t)}\langle O(t)k, \phi \rangle + \int_{S(t)}\langle
 O(t^{2})\dot k, \phi \rangle .
\end{equation}
Here $P_{2}(\phi)$ is obtained by integrating the $D_{x}^{2}$ term on
 the right in 
(4.17) by parts over $S(t)$. Thus $P_{2}(\phi)$, and more generally,
 $P_{k}(\phi)$ 
denote differential operators of order $k$ on $\phi$ with coefficients
 depending 
on $g$ and $g_{1}$ and their derivatives up to order 2 and so at least
 continuous 
up to $\bar \partial \Omega$. We use these expressions generically, so
 their 
exact form may change from line-to-line below. Note also there are no
 boundary 
terms at $\partial S(t)$ arising from the integration by parts, by the
 vanishing 
hypothesis on $\partial C$. 

    For the terms on the right in (4.24) one then has
$$\int_{S(t)}\langle O(t)k, \phi \rangle = t\int_{S(t)}\langle k, 
P_{0}(\phi) \rangle ,$$
while, since $A = O(t)$ and $H = O(t)$, 
$$\int_{S(t)}\langle O(t^{2})\dot k, \phi \rangle = 
t^{2}\int_{S(t)}\langle \dot k, P_{0}(\phi) \rangle 
= t^{2}\frac{d}{dt}\int_{S(t)}\langle k, P_{0}(\phi) \rangle - 
t^{2}\int_{S(t)}\langle k, P_{1}(\phi) \rangle .$$
Similarly, for the terms on the left in (4.24), one has
$$\int_{S(t)}\langle \dot k, \phi \rangle  = 
\frac{d}{dt}\int_{S(t)}\langle k, \phi \rangle  - t \int_{S(t)}\langle
 k, P_{1}(\phi) \rangle ,$$
while
$$\int_{S(t)}\langle \ddot k, \phi \rangle = 
\frac{d^{2}}{dt^{2}}\int_{S(t)}\langle k, \phi \rangle 
- 2t\frac{d}{dt}\int_{S(t)}\langle k, P_{1}(\phi) \rangle + 
\int_{S(t)}\langle k, P_{1}(\phi) \rangle  + t\int_{S(t)}\langle k,
 P_{2}(\phi) \rangle .$$

   Now let 
$$f = f(t) = \int_{S(t)}\langle k, \phi \rangle  .$$
Then the computations above give
\begin{equation} \label{e4.25}
t\ddot f - (n-1)\dot f = t\int_{S(t)}\langle k, P_{2}(\phi) \rangle +
 (1 + t^{2})\int_{S(t)}
\langle  k, P_{1}(\phi) \rangle
\end{equation}
$$+ \frac{d}{dt}\int_{S(t)}t^{2}\langle k, P_{0}(\phi) \rangle 
+ \frac{d}{dt}\int_{S(t)}t\langle k, P_{1}(\phi) \rangle .$$

   First observe that 
\begin{equation}\label{e4.26}
\int_{S(t)}\langle k, \phi \rangle  = o(t^{n}) \Rightarrow
 \int_{S(t)}\langle k, 
P_{k}(\phi) \rangle = o(t^{n}),
\end{equation}
for all $C^{\infty}$ forms $\phi$ vanishing to infinite order at
 $\partial C$. 
For if the left side of (4.26) holds, then $\int_{S(t)}\langle k,
 \partial^{k}\phi 
\rangle  = o(t^{n})$, since the hypotheses on $\phi$ are closed under
 differentiation. 
The coefficients of $P_{k}$ are at least continuous, and it is
 elementary to verify 
that if $\int_{S(t)}\langle k, \partial^{k}\phi \rangle = o(t^{n})$,
 then 
$\int_{S(t)}\langle k, \phi \partial^{k}\phi \rangle = o(t^{n})$, for
 any function 
$\phi$ continuous on $\bar C$. Note that the same result holds
 with $p$ in 
place of $n$, for any $p < \infty$. 

  It follows from (4.26) and the initial hypothesis (4.22) that the
 first two terms 
on the right in (4.25) are $o(t^{n})$ as $t \rightarrow 0$. Since
 $t\ddot f - 
(n-1)\dot f = t^{n}\frac{d}{dt}(\frac{\dot f}{t^{n-1}})$, this gives
$$\frac{d}{dt}(\frac{\dot f}{t^{n-1}}) = o(1) +
 t^{-n}\frac{d}{dt}\int_{S(t)}t\langle k, 
P_{1}(\phi) \rangle + t^{-n}\frac{d}{dt}\int_{S(t)}t^{2}\langle k,
 P_{0}(\phi) \rangle .$$
Integrating from $0$ to $t$ implies
$$\frac{\dot f}{t^{n-1}} = o(t) + t^{-n+1}\int_{S(t)}\langle k,
 P_{1}(\phi) \rangle 
+ n\int_{0}^{t}t^{-n}\int_{S(t)}\langle k, P_{1}(\phi) \rangle + c_{1}
 = o(t) 
+ c_{1} ,$$
where $c_{1}$ is a constant. A further integration using (4.26) again
 gives
\begin{equation} \label{e4.27}
f = o(t^{n+1}) + c_{1}'t^{n} + c_{2},
\end{equation}
where $c_{1}' = \frac{c_{1}}{n}$. Once more by (4.22), this implies
 that
$$f = o(t^{n+1}).$$
Note the special role played by the indicial root $n$ here; if instead
 one had only 
$k = O(t^{n})$, then the argument above does not give $k = O(t^{n+1})$
 weakly. 

  This first estimate holds in fact for any given trace-free $\phi$
 which is $C^{2}$ 
on $\bar C$, and vanishing to first order on $\partial C$.
 Working in the 
same way with the trace equation (4.20) shows that the same result
 holds for pure 
trace terms. In particular, it follows that
\begin{equation} \label{e4.28}
k = o(t^{n+1}) \ {\rm weakly} .
\end{equation}

  One now just repeats this argument inductively, with the improved
 estimate 
(4.28) in place of (4.22), using (4.26) inductively. Note that each
 inductive 
step requires higher differentiability of the test function $\phi$ and
 its higher 
order vanishing at $\partial C$. 
{\endproof}

  Lemma 4.2 proves that $k = k_{\alpha\beta} = O(t^{\nu})$ weakly, for
 any 
$\nu < \infty$. As discussed in \S 3, the transition from geodesic
 boundary 
coordinates to $H$-harmonic coordinates is $C^{2,\alpha}$ and hence
\begin{equation}\label{e4.29}
h = h_{\alpha\beta} = O(t^{\nu}),
\end{equation}
weakly, with the level sets $S(t)$ replaced by $S_{\tau}$. Next,
 as in 
Remark 3.5 and the proof of Theorem 3.1, the equations (4.11)-(4.13)
 satisfy 
elliptic estimates, and elliptic regularity in weighted H\"older
 spaces, cf.~[26], [20], 
shows that the weak decay (4.29) implies strong or pointwise decay,
 i.e.~(4.14) holds. 
The proof of Theorem 4.1 and thus Theorem 1.2 is now completed as
 before in the 
$C^{\infty}$ smooth case. 
{\endproof}

\begin{remark} \label{r 4.3}
{\rm In [3, Thm.~3.2], a proof of unique continuation of conformally
 compact 
Einstein metrics was given in dimension 4, using the fact that the 
compactified metric $\widetilde g$ in (1.1) satisfies the Bach
 equation, 
together with the Calder\'on uniqueness theorem. However, the proof 
in [3] used harmonic coordinates; as discussed in \S 2, such
 coordinates 
do not preserve the Cauchy data. The first author is grateful to 
Robin Graham for pointing this out. Theorem 1.2 thus corrects this 
error, and generalizes the result to any dimension. }
\end{remark}

  For the work to follow in \S 5, we note that Theorem 4.1 also holds
 for 
linearizations of the Einstein equations, i.e.~forms $k$ satisfying 
\begin{equation}\label{e4.30}
\frac{d}{dt}(\Ric_{g+tk} + n(g+tk))|_{t=0} = 0. 
\end{equation}
Thus, if $k$ satisfies (4.30) and the analog of (4.7), i.e.~$|k| =
 o(t^{n})$, then 
$k$ is pure gauge in $\Omega$, in that $k = \delta^{*}Z$, where $Z$ is
 a vector 
field on $\Omega$ with $Z = 0$ on $\partial \Omega$. The proof of this
 is exactly 
the same as the proof of Theorem 4.1, replacing the finite difference
 $k = g_{1} - g_{0}$ by an infinitesimal difference. 

  This has the following consequence:
\begin{corollary}\label{c4.4}
Let $(M, g)$ be a conformally compact Einstein manifold with metric $g$
 having a 
$C^{3,\alpha}$ geodesic compactification. Suppose the topological
 condition 
{\rm (1.5)} holds, i.e.~$\pi_{1}(M, \partial M) = 0$. 

  If $k$ is an infinitesimal Einstein deformation on $M$ as in {\rm
 (4.30)}, 
in divergence-free gauge, i.e.
\begin{equation}\label{e4.31}
\delta k = 0,
\end{equation}
with $k = o(t^{n})$ on approach to $\partial M$, then
$$k = 0 \ \ {\rm on} \ \ M.$$
\end{corollary}

{\bf Proof:}  The topological condition (1.5), together with the same 
analytic continuation argument at the end of the proof of Theorem 3.1, 
implies that $k$ is pure gauge globally on $M$, in that $k = \delta^{*}Z$ 
on $M$ with $Z = 0$ on $\partial M$. (Recall that (1.5) implies that 
$\partial M$ is connected). From (4.31), one then has
$$\delta \delta^{*}Z = 0,$$
on $M$. Pairing this with $Z$ and integrating over $B(t)$, it follows
 that 
$$\int_{B(t)}|\delta^{*}Z|^{2} = \int_{S(t)}\delta^{*}Z(Z, N),$$
where $N$ is the unit outward normal. Since $|Z|_{g}$ is bounded and 
$|\delta^{*}Z|vol(S(t)) = o(1)$, (since $|k| = o(t^{n})$), it follows
 that 
$$\int_{M}|\delta^{*}Z|^{2} = 0,$$
which gives the result.
{\endproof}

  Of course, analogs of these results also hold for bounded domains,
 via the 
proof of Theorem 3.1; the verification is left to the reader. 

\begin{remark} \label{r4.5}
{\rm The analogue of Proposition 3.7 most likely also holds in the
 setting of 
conformally compact metrics, for fields $\tau$ whose Euler-Lagrange 
equation is a diagonal system of Laplace-type operators to leading 
order, as in \eqref{e3.46} or \eqref{e3.47}. The proof of this is 
basically the same as that of Proposition 3.7, using the proof of 
Theorem 1.2 and with the Mazzeo unique continuation result in place 
of that of Calder\'on. However, we will not carry out the details 
of the proof here. }
\end{remark}

\section{Isometry Extension and the Constraint Equations.}
\setcounter{equation}{0}

  In this section, we prove Theorem 1.3 that continuous groups of
 isometries at the 
boundary extend to isometries in the interior of complete conformally
 compact Einstein 
metrics and relate this issue in general to the constraint equations
 induced by the 
Gauss-Codazzi equations.  

  We begin with the following elementary consequence of Theorem 4.1. 

\begin{proposition} \label{p5.1}
Let $(\Omega, g)$ be a $C^{n}$ polyhomogeneous conformally compact
 Einstein metric 
on a domain $\Omega \simeq B^{n+1}$ with boundary metric $\gamma$ on
 $\partial 
\Omega \simeq B^{n}$. Suppose $X$ is a Killing field on $(\partial
 \Omega, \gamma)$ 
and 
\begin{equation} \label{e5.1}
{\mathcal L}_{X}g_{(n)} = 0,
\end{equation}
where $g_{(n)}$ is the $n^{\rm th}$ term in the Fefferman-Graham
 expansion {\rm (4.3)} 
or {\rm (4.4)}.

  Then $X$ extends to a Killing field on $(\Omega, g)$. 
\end{proposition}

{\bf Proof:} 
Extend $X$ to a smooth vector field on $\Omega$ by requiring $[X, N] =
 0$, where 
$N = \nabla \log t$ and $t$ is the geodesic defining function
 determined by $g$ and 
$\gamma$. Let $\phi_{s}$ be the corresponding 1-parameter group of
 diffeomorphisms 
and set $g_{s} = \phi^{*}_{s}g$. Then $t$ is the geodesic defining
 function for 
$g_{s}$ for any $s$, and the pair $(g, g_{s})$ satisfy the hypotheses
 of Theorem 4.1. 
Theorem 4.1 then implies that $g_{s}$ is isometric to $g$, i.e.~there 
exist diffeomorphisms $\psi_{s}$ of $\Omega$, equal to the identity on 
$\partial \Omega$, such that $\psi_{s}^{*}\phi_{s}^{*}g = g$. Thus 
$\phi_{s}\circ \psi_{s}$ is a 1-parameter group of isometries of $g$
 defined in 
$\Omega$, with $Y$ the corresponding Killing field. (In fact, $Y = X$,
 since any Killing 
field $Y$ tangent to $\partial\Omega$ preserves the geodesics tangent
 to $N$, and so 
$[Y, N] = 0$. This determines $Y$ uniquely in terms of its value at
 $\partial \Omega$. 
Since $X$ satisfies the same equation with the same initial value, this
 gives the claim). 

{\endproof}

  We point out that the the same result, and proof, also hold in the
 case of 
Einstein metrics on bounded domains, via Theorem 3.1; the condition
 (5.1) is 
of course replaced by ${\mathcal L}_{X}A = 0$. For some examples and
 discussion 
in the bounded domain case, see [1], [2]. 

  Suppose now that $(M, g)$ is a (global) conformally compact Einstein
 metric and 
there is a domain $\Omega$ as in Proposition 5.1 contained in $M$ on
 which (5.1) 
holds. Then by analytic continuation as discussed at the end of the
 proof of Theorem 
3.1, $X$ extends to a local Killing field on all of $M$, i.e.~$X$
 extends to a Killing 
field on the universal cover $\widetilde M$. In particular, if the
 condition \eqref{e1.5} holds, i.e.~$\pi_{1}(M, \partial M) = 0$, 
then $X$ extends to a global Killing field on $M$. Again, the same 
result holds in the context of bounded domains. 

\begin{remark}\label{r5.2}
{\rm A natural analogue of Proposition 5.1 holds for conformal Killing
 fields on 
$(\partial \Omega, \gamma)$, i.e.~vector fields which preserve the
 conformal class 
$[\gamma]$ at conformal infinity. Such vector fields satisfy the
 conformal Killing 
equation
\begin{equation}\label{e5.2}
\hat{\mathcal L}_{X}\gamma = {\mathcal L}_{X}\gamma - \frac{tr(\mathcal
 L_{X}\gamma)}{n}
\gamma = 0 .
\end{equation}
Namely, since we are working locally, it is well-known - and easy to
 prove - that any 
non-vanishing conformal Killing field is Killing with respect to a
 conformally related 
metric $\widetilde \gamma = \lambda^{2}\gamma$, so that 
$${\mathcal L}_{X}\widetilde \gamma = 0 .$$
Hence, if ${\mathcal L}_{X}\widetilde g_{(n)} = 0$, then Proposition
 5.1 implies that 
$X$ extends to a Killing field on $\Omega$. 

   One may express $\widetilde g_{(n)}$ in terms of $\lambda$ and the
 lower order 
terms $g_{(k)}$, $k < n$ in the Fefferman-Graham expansion (4.3)-(4.4);
 however, the 
expressions become very complicated for $n$ even and large, cf.~[17].
 Thus, while the 
equation (5.2) is conformally invariant, the corresponding conformally
 invariant 
equation for $g_{(n)}$ will be complicated in general. }
\end{remark}

   Next we consider the constraint equations (4.5) in detail, i.e.
\begin{equation}\label{e5.3}
\delta \tau_{(n)} = 0 \ \ {\rm and} \ \ tr \,\tau_{(n)} = a_{(n)},
\end{equation}
where $\tau_{(n)} = g_{(n)} + r_{(n)}$; $r_{(n)}$ and $a_{(n)}$ are
 explicitly determined 
by the boundary metric $\gamma = g_{(0)}$ and its derivatives up to
 order $n$. Both vanish 
when $n$ is odd. 

  As will be seen below, the most important issue is the divergence
 constraint in 
(5.3), which arises from the Gauss-Codazzi equations. To see this, in
 the setting 
of \S 4, on $S(t) \subset (M, g)$, the Gauss-Codazzi equations are
\begin{equation} \label{e5.4}
\delta(A - Hg) = -\Ric(N, \cdot),
\end{equation}
as 1-forms on $S(t)$; here $N = -t\partial_{t}$ is the unit outward
 normal. The same 
equation holds on a geodesic compactification $(M, \bar g)$. If $g$ is
 Einstein, then 
$\Ric(N, \cdot) = \bar \Ric(\bar N, \cdot) = 0$; the latter equality
 follows from (4.9). 
The equation (5.4) holds for all $t$ small, and differentiating $(n-1)$
 times with 
respect to $t$ gives rise to the divergence constraint in (5.3). 

   The Gauss-Codazzi equations are not used in the derivation and
 properties of 
the Fefferman-Graham expansion (4.3)-(4.4) per se. The derivation of
 these equations 
involves only the tangential $(ij)$ part of the Ricci curvature. The
 asymptotic 
behavior of the normal $(00)$ part of the Ricci curvature gives rise to
 the trace 
constraint in (5.3), cf.~(4.19)-(4.20). 

\medskip

  Let ${\mathcal T}$ be the space of pairs $(g_{(0)}, \tau_{(n)})$
 satisfying (5.3). 
If $\tau_{(n)}^{0}$ is any fixed solution of (5.3), then any other
 solution with 
the same $g_{(0)}$ is of the form $\tau_{(n)} = \tau_{(n)}^{0} + \tau$,
 where 
$\tau$ is transverse-traceless with respect to $g_{(0)}$. (Of
 course if $n$ 
is odd, one may take $\tau_{(n)}^{0} = 0$). The space ${\mathcal T}$ 
 naturally projects onto $\Met(\partial M)$ with fiber at $\gamma$ an 
affine space of symmetric tensors 
and is a subset of the product 
$\Met(\partial M) \times {\mathbb S}^{2}(\partial M) \simeq
 T(\Met(\partial M))$. Let 
\begin{equation}\label{e5.5}
\pi: {\mathcal T} \rightarrow \Met(\partial M)
\end{equation}
be the projection onto the base space $\Met(\partial M)$, (the first
 factor projection). 

  By the discussion in \S 4, $(g_{(0)}, \tau_{(n)}) \in  {\mathcal T}$
 if and only 
if the corresponding pair $(g_{(0)}, g_{(n)})$ determine a formal
 polyhomogenous solution 
to the Einstein equations near conformal infinity, i.e.~formal series
 solutions 
containing $\log$ terms, as in (4.3)-(4.4). In fact, if  $g_{(0)}$ and
 $g_{(n)}$ are 
real-analytic on $\partial M$, a result of Kichenassamy [24] implies
 that the series 
(4.3) or (4.4) converges, and gives an Einstein metric $g$, defined in
 a neighborhood 
of $\partial M$. The metric $g$ is complete near $\partial M$ and has a
 conformal 
compactification inducing the given data $(g_{(0)}, g_{(n)})$ on
 $\partial M$. 
Here we recall from the discussion in \S 4 that all coefficients of the
 expansion 
(4.3) or (4.4) are determined by $g_{(0)}$ and $g_{(n)}$. 

\medskip

In this regard, consider the following:

  {\bf  Problem.} Is $\pi: {\mathcal T} \rightarrow \Met(\partial M)$ an
 open map? 
Thus, given any $(g_{(0)}, \tau_{(n)}) \in {\mathcal T}$ and any
 boundary metric 
$\widetilde g_{(0)}$ sufficiently close to $g_{(0)}$, does there exist 
$\widetilde \tau_{(n)}$ close to $\tau_{(n)}$ such that $(\widetilde
 g_{(0)}, 
\widetilde \tau_{(n)}) \in {\mathcal T}$.

  Although $\pi$ is obviously globally surjective, the problem above is
 whether 
$\pi$ is locally surjective. For example, a simple fold map $x
 \rightarrow x^{3}-x$ 
is not locally surjective near $\pm\sqrt{3}/3$. Observe that the trace
 condition in 
(5.3) imposes no constraint on $g_{(0)}$; given any $g_{(0)}$, it is
 easy to find 
$g_{(n)}$ such that $tr_{g_{(0)}}(g_{(n)}+r_{(n)}) = a_{(n)}$; this
 equation can readily 
be solved algebraically for many $g_{(n)}$. 

  By the inverse function theorem, it suffices, (and is probably also
 necessary), to 
examine the problem above at the linearized level. However the
 linearization of the 
divergence condition in (5.3) gives a non-trivial constraint on the
 variation $h_{(0)}$ 
of $g_{(0)}$. Namely, the linearization in this case gives 
\begin{equation} \label{e5.6} 
\delta'(\tau_{(n)}) + \delta(\tau_{(n)})' = 0 ,
\end{equation}
where $\delta ' = \frac{d}{du}\delta_{g_{(0)}+uh_{(0)}}$, and similarly
 for 
$(\tau_{(n)})'$.

  Whether (5.6) is solvable for any $h_{(0)} \in S^{2}(\partial M)$ 
depends on the data $g_{(0)}$ and $g_{(n)}$. For example, it is
 trivially solvable when 
$\tau_{(n)} = 0$. For compact $\partial M$, one has 
\begin{equation} \label{e5.7}
\Omega^{1}(\partial M) = Im \delta \oplus Ker \delta^{*}, 
\end{equation}
where $\Omega^{1}$ is the space of 1-forms, so that solvability in
 general requires that 
\begin{equation}\label{e5.8}
\delta'(\tau_{(n)}) \in Im \delta = (Ker \delta^{*})^{\perp}.
\end{equation}
Of course $Ker \delta^{*}$ is exactly the space of Killing fields on
 $(\partial M, \gamma)$, 
and so this space serves as a potential obstruction space.  

  Clearly then $\pi$ is locally surjective when $(\partial M, g_{(0)})$
 has no Killing fields. 
On the other hand, it is easy to construct examples where $(\partial M,
 \gamma)$ does 
have Killing fields and $\pi$ is not locally surjective: 

\begin{example} \label{ex5.3}
{\rm  Let $(\partial M, g_{(0)})$ be the flat metric on the $n$-torus
 $T^{n}$, $n \geq 3$, 
and define $g_{(n)} = -(n-2)(d\theta^{2})^{2} +(d\theta^{3})^{2} +
 \cdots + (d\theta^{n})^{2}$. 
Then $g_{(n)}$ is transverse-traceless with respect to $g_{(0)}$. Let
 $f = f(\theta^{1})$. 
Then $\hat g_{(n)} = fg_{(n)}$ is still transverse-traceless with 
respect to $g_{(0)}$, so that 
$(g_{(0)}, \hat g_{(n)}) \in {\mathcal T}$, at least for $n$ odd. 

   It is then not difficult to see via a direct calculation, or more
 easily via 
Proposition 5.4 below, that (5.8) does not hold, so that $\pi$ is not
 locally surjective. }
\end{example}

   Next we relate these two issues, i.e.~the general solvability of the
 divergence 
constraint (5.6) and the extension of Killing fields on the boundary
 into the bulk. 
The following result holds for general $\phi \in S^{2}(\partial M)$
 with $\delta \phi = 0$ on $\partial M$. 

\begin{proposition} \label{p5.4}
If $X$ is a Killing field on $(\partial M, \gamma)$, with $\partial M$
 compact, then 
\begin{equation}\label{e5.9}
\int_{\partial M}\langle {\mathcal L}_{X}\tau_{(n)}, h_{(0)} \rangle
 dV= 
-2\int_{\partial M}\langle \delta'(\tau_{(n)}), X \rangle dV,
\end{equation}
where $\delta' = \frac{d}{ds}\delta_{\gamma + sh_{(0)}}$. 
In particular, {\rm (5.1)} holds for all Killing fields on $(\partial
 M, \gamma)$ if 
and only if the linearized divergence constraint vanishes, i.e.~{\rm
 (5.6)} holds for all $h$. 
\end{proposition}

{\bf Proof:}
Since $X$ is a Killing field on $(\partial M, \gamma)$, one has
\begin{equation} \label{e5.10}
\int_{\partial M}\langle {\mathcal L}_{X}\tau, h \rangle dV_{\gamma} =
 - 
\int_{\partial M}\langle \tau, {\mathcal L}_{X}h \rangle dV_{\gamma}.
\end{equation}
Setting $\gamma_{s} = \gamma + sh$, the divergence theorem gives
\begin{equation}\label{e5.11}
0 = \int_{\partial M}\delta_{\gamma_{s}}(\tau(X))dV_{\gamma_{s}} = 
\int_{\partial M}\langle \delta_{\gamma_{s}}\tau, X \rangle
 dV_{\gamma_{s}} 
- {\tfrac{1}{2}}\int_{\partial M}\langle \tau, {\mathcal
 L}_{X}\gamma_{s} 
\rangle dV_{\gamma_{s}},
\end{equation}
where the second equality is a simple computation from the definitions;
 the 
inner products are with respect to $\gamma_{s}$. Taking the derivative
 with 
respect to $s$ at $s = 0$, and using the facts that $X$ is Killing and 
$\delta(\tau) = 0$, it follows that
\begin{equation}\label{e5.12}
\int_{\partial M}\langle \delta' \tau, X \rangle dV 
- {\tfrac{1}{2}}\int_{\partial M}\langle \tau, {\mathcal L}_{X}h 
\rangle dV = 0.
\end{equation}
Combining this with (5.10) then gives (5.9); note that 
${\mathcal L}_{X}r_{(n)} = 0$ in this case, since $r_{(n)}$ 
is determined by the boundary metric. 

  To prove the last statement, by (5.9), (5.1) holds if and only if 
$\int_{\partial M}\langle \delta'(\tau_{(n)}), X \rangle = 0$, for all 
variations $h$. If (5.6) holds, then $\delta'(\tau_{(n)}) = \delta
 h_{(n)}'$, 
for some $h_{(n)}'$ and so $\int_{\partial M}\langle
 \delta'(\tau_{(n)}), X \rangle 
= \int_{\partial M}\langle h_{(n)}', \delta^{*}X \rangle = 0$, since
 $X$ is 
Killing. The converse of this argument holds equally well. 

{\endproof}

  Proposition 5.4 implies that in general, Killing fields on $\partial
 M$ do not 
extend to Killing fields in a neighborhood of $\partial M$,
 (cf.~Example 5.3). 
(Exactly the same result and proof hold in the bounded domain case,
 when the term 
$\tau_{(n)}$ is replaced by $A - Hg$). 

  Now as noted above, whether isometry extension holds or not depends
 on the 
term $\tau_{(n)} = g_{(n)}+r_{(n)}$, or more precisely on the relation
 of the 
boundary metric $g_{(0)}$ with $\tau_{(n)}$. For Einstein metrics which
 are 
globally conformally compact, the term $\tau_{(n)}$ is determined, up
 to a finite 
dimensional moduli space, by the boundary metric $g_{(0)}$; (this is
 discussed 
further below). Thus, whether isometry extension holds or not is quite
 a delicate 
issue; if so, it must depend crucially on the global structure of $(M,
 g)$. 

\medskip

  Before beginning the proof of Theorem 1.3, we first need to discuss
 some 
background material from [5]-[6].

  Let $E_{AH}$ be the space of conformally compact, or equivalently
 asymptotically 
hyperbolic Einstein metrics on $M$ which have a $C^{\infty}$
 polyhomogeneous
conformal compactification with respect to a fixed smooth defining
 function $\rho$, 
as in (1.1). In [5], it is shown that $E_{AH}$ is a smooth, infinite
 dimensional 
manifold. One has a natural smooth boundary map 
\begin{equation}\label{e5.13}
\Pi: E_{AH} \rightarrow \Met(\partial M),
\end{equation}
sending $g$ to its boundary metric $\gamma$. 

  The moduli space ${\mathcal E}_{AH}$ is the quotient 
$E_{AH}/{\mathcal D}_{1}$, where ${\mathcal D}_{1}$ is the 
group of smooth (polyhomogeneous) diffeomorphisms $\phi$ of $M$ 
equal to the identity on $\partial M$. Thus, $g' \sim g$ if 
$g' = \phi^{*}g$, with $\phi \in {\mathcal D}_{1}$. Changing the 
defining function $\rho$ in (1.1) changes the boundary metric conformally. 
Also, if $\phi \in {\mathcal D}_{1}$ then $\rho \circ \phi$ is another 
defining function, and all defining functions are of this form near 
$\partial M$. Hence if ${\mathcal C}$ denotes the space of smooth 
conformal classes of metrics on $\partial M$, then the boundary map 
(5.13) descends to a smooth map 
\begin{equation}\label{e5.14}
\Pi: {\mathcal E}_{AH} \rightarrow {\mathcal C}
\end{equation}
independent of the defining function $\rho$. The boundary map $\Pi$
 in (5.14) is Fredholm, of Fredholm index 0. 

 The linearization of the Einstein operator $\Ric_{g} + ng$ at an
 Einstein metric 
$g$ is given by 
\begin{equation}\label{e5.15}
\hat L = (\Ric_{g} + ng)' = \tfrac{1}{2}D^{*}D -  R - \delta^{*}\beta, 
\end{equation}
acting on the space of symmetric 2-tensors $S^{2}(M)$ on $M$, cf.~[10].
 Here, (as 
in \S 3), $\beta$ is the Bianchi operator, $\beta(h) = \delta h +
 \frac{1}{2}d tr h$, 
Thus, $h \in T_{g}E_{AH}$ if and only if 
$$\hat L(h) = 0.$$
The operator $\hat L$ is not elliptic, due to the $\delta^{*}\beta$
 term. As is 
well-known, this arises from the diffeomorphism group, and to obtain an
 elliptic 
linearization, one needs a gauge choice to break the diffeomorphism
 invariance 
of the Einstein equations. We will use a slight modification of the
 Bianchi 
gauge introduced in [11]. 

  To describe this, given any fixed $g_{0} \in E_{AH}$ with geodesic
 defining 
function $t$ and boundary metric $\gamma_{0}$, let $\gamma$ be a
 boundary metric 
near $\gamma_{0}$ and define the hyperbolic cone metric $g_{\gamma}$ on
 $\gamma$ 
by setting 
$$g_{\gamma} = t^{-2}(dt^{2} + \gamma);$$
$g_{\gamma}$ is defined in a neighborhood of $\partial M$. Next, set
\begin{equation}\label{e5.16}
g(\gamma) = g_{0} + \eta(g_{\gamma} - g_{\gamma_{0}}),
\end{equation}
where $\eta$ is a non-negative cutoff function supported near $\partial
 M$ with 
$\eta = 1$ in a small neighborhood of $\partial M$. Any conformally
 compact metric 
$g$ near $g_{0}$, with boundary metric $\gamma$ then has the form 
\begin{equation}\label{e5.17}
g = g(\gamma) + h,
\end{equation}
where $|h|_{g_{0}} = O(t^{2})$; equivalently $\bar h = t^{2}h$
 satisfies $\bar h_{ij} 
= O(t^{2})$ in any smooth coordinate chart near $\partial M$. The space
 of such 
symmetric bilinear forms $h$ is denoted by ${\mathbb S}_{2}(M)$ and the
 space of 
metrics $g$ of the form (5.17) is denoted by $\Met_{AH}$. 

  The Bianchi-gauged Einstein operator, (with background metric
 $g_{0}$), is defined 
by
\begin{equation}\label{e5.18}
\Phi_{g_{0}}: \Met_{AH} \rightarrow {\mathbb S}_{2}(M)
\end{equation}
$$\Phi_{g_{0}}(g) = \Phi (g(\gamma) + h) = \Ric_{g} + ng + 
(\delta_{g})^{*}\beta_{g(\gamma)}(g),$$
where $\beta_{g(\gamma)}$ is the Bianchi operator with respect to
 $g(\gamma)$. By 
[11, Lemma I.1.4],
\begin{equation}\label{e5.19}
Z_{AH} \equiv  \Phi^{-1}(0)\cap\{\Ric <  0\} \subset E_{AH}, 
\end{equation}
where $\{\Ric < 0\}$ is the open set of metrics with negative Ricci 
curvature. In fact, if $g \in E_{AH}$ is close to $g_{0}$, and $\Phi (g) 
= 0$, then $\beta_{g(\gamma)}(g) = 0$ and moreover
\begin{equation}\label{e5.20}
\delta_{g(\gamma)}(g) = 0 \ \ {\rm and} \ \ tr_{g(\gamma)}(g) = 0.
\end{equation}
The space $Z_{AH}$ is a local slice for the action of ${\mathcal D}_{1}$ 
on $E_{AH}$: for any $g\in E_{AH}$ near $g_{0}$, there exists a 
diffeomorphism $\phi \in {\mathcal D}_{1}$ such that $\phi^{*}g\in 
Z_{AH}$, cf.~again [11].

  The linearization of $\Phi$ at $g_{0} \in E_{AH}$ with respect to
 the $2^{\rm nd}$ 
variable $h$ has the simple form
\begin{equation}\label{e5.21}
(D_{2}\Phi)_{g_{0}}(\dot h) = \tfrac{1}{2}D^{*}D \dot h -
  R_{g_{0}}(\dot h), 
\end{equation}
while the variation of $\Phi$ at $g_{0}$ with respect to the $1^{\rm
 st}$ variable 
$g(\gamma)$ has the form 
\begin{equation}\label{e5.22}
(D_{1}\Phi)_{g_{0}}(\dot g(\gamma)) = (D_{2}\Phi)_{g_{0}}(\dot
 g(\gamma)) 
- \delta_{g_{0}}^{*}\beta_{g_{0}}(\dot g(\gamma)) = (\Ric_{g} +
 ng)'(\dot g(\gamma)),
\end{equation}
as in (5.15). Clearly $\dot g(\gamma) = \eta t^{-2}\dot \gamma$. 
  The kernel of the elliptic self-adjoint linear operator
\begin{equation}\label{e5.23} 
L = \tfrac{1}{2}D^{*}D -  R 
\end{equation}
acting on the $2^{\rm nd}$ variable $h$, represents the space of
 non-trivial 
infinitesimal Einstein deformations vanishing on $\partial M$. Let $K$
 denote the 
$L^{2}$ kernel of $L$. This is the same as the kernel of $L$ on
 ${\mathbb S}_{2}(M)$, 
cf.~[11], [26]. An Einstein metric $g_{0} \in E_{AH}$ is called {\it
 non-degenerate} 
if 
\begin{equation} \label{e5.24}
K = 0. 
\end{equation}

  For $g_{0} \in {\mathcal E}_{AH}$ the kernel $K = K_{g_{0}}$ equals
 the kernel of the linear map $D\Pi: T_{g_{0}}{\mathcal E}_{AH} \rightarrow 
T_{\Pi(g_{0})}{\mathcal C}$. Hence, $g_{0}$ is non-degenerate if and 
only if $g_{0}$ is a regular point of the boundary map $\Pi$ in which 
case $\Pi$ is a local diffeomorphism near $g_{0}$. From now on, we 
denote $g_{0}$ by $g$. 

  By the regularity result of Chru\'sciel et al.~[16], any $\kappa \in K$ 
has a $C^{\infty}$ smooth polyhomogeneous expansion, analogous to the 
Fefferman-Graham expansion (4.3)-(4.4), with leading order terms satisfying 
\begin{equation}\label{e5.25}
\kappa = O(t^{n}), \ \ \kappa(N,Y) = O(t^{n+1}), \ \ \kappa (N,N) = 
O(t^{n+1+\mu}),
\end{equation}
where $N = -t\partial_{t}$ is the unit outward normal vector to the 
$t$-level set $S(t)$, $Y$ is any $g$-unit vector tangent to $S(t)$ 
and $\mu > 0$; cf.~also [28, Prop.~5]. Here $\kappa = O(t^{n})$ means 
$|\kappa|_{g} = O(t^{n})$. Also by an argument similar to the one
leading to (5.20), any $\kappa \in K$ is 
transverse-traceless, i.e.
\begin{equation}\label{e5.26}
\delta \kappa = tr \kappa = 0.
\end{equation}

  Given this background, we are now ready to begin the proof of Theorem
 1.3. 

\medskip

{\bf Proof of Theorem 1.3.}

\medskip
 
  Let $\bar g = t^{2}g$ be a geodesic compactification of $g$ with 
boundary metric $\gamma$. By the boundary regularity result of [16], 
$\bar g$ is $C^{\infty}$ polyhomogeneous on $\bar M$. It suffices to
 prove 
Theorem 1.3 for arbitrary 1-parameter subgroups of the isometry group
 of 
$(\partial M, \gamma)$. Thus, let $\phi_{s}$ be a local 1-parameter
 group 
of isometries of $\gamma$ with $\phi_{0} = id$, so that
$$\phi_{s}^{*}\gamma  = \gamma .$$ 
The diffeomorphisms $\phi_{s}$ of $\partial M$ may be extended to 
diffeomorphisms of $M$, so that the curve  
\begin{equation} \label{e5.27}
g_{s} = \phi_{s}^{*}g
\end{equation}
is a smooth curve in $E_{AH}$. By construction then, $\Pi[g_{s}] =
 [\gamma]$, 
so that $[h] = [\frac{dg_{s}}{ds}] \in Ker D\Pi$, for $\Pi$ as in
 (5.14). One 
may then alter the diffeomorphisms $\phi_{s}$ by composition with
 diffeomorphisms 
in ${\mathcal D}_{1}$ if necessary, so that $h = \frac{dg_{s}}{ds} \in
 K_{g}$, 
where $K_{g}$ is the kernel in (5.24). Denoting $h = \kappa$, it
 follows that 
\begin{equation}\label{e5.28}
\kappa = \delta^{*}X,
\end{equation}
where $X = d\phi_{s}/ds$ is smooth up to $\bar M$. 

  Thus it suffices to prove that $\delta^{*}X = 0$, since this will
 imply that 
$g_{s} = g$, (when $g_{s}$ is modified by the action of ${\mathcal
 D}_{1}$). If 
$K_{g} = 0$, i.e.~if $g$ is a regular point of the boundary map $\Pi$,
 then 
this is now obvious, (from the above), and proves the result in this
 special 
case; (the proof in this case requires only that $(M, g)$ be
 $C^{2,\alpha}$ 
conformally compact).

  We give two different, (although related), proofs of Theorem 1.3, one
 conceptual and 
one more computational. The first, conceptual, proof involves an
 understanding of the 
cokernel of the map $D\Pi_{g}$ in $\Met(\partial M)$, and so one first
 needs to give an 
explicit description of this cokernel. To begin, recall the derivative
\begin{equation} \label{e5.29}
(D\Phi)_{g}: T_{g}\Met_{AH}(M) \rightarrow  T_{\Phi(g)}{\mathbb
 S}_{2}(M). 
\end{equation}
Via (5.17), one has $T_{g}\Met_{AH} = T_{\gamma}\Met(\partial M)\oplus
 T_{h}{\mathbb S}_{2}(M)$ 
and the derivative with respect to the second factor is given by
 (5.21). If $K = 0$, then 
$D_{2}\Phi$ is surjective at $g$, (since $D_{2}\Phi$ has index 0, and we 
recall that the kernel and cokernel here are equal to their $L^2$ counterparts), 
and hence so is $D\Pi$. 
In general, to understand $Coker D\Pi$, we show that $D\Phi$ is always
 surjective; this 
follows from the claim that for any non-zero $\kappa\in K$ there is a
 tangent vector 
$\dot g(\gamma) \in T_{\gamma}\Met(\partial M) \subset T_{g}\Met_{AH}$
 such that
\begin{equation} \label{e5.30}
\int_{M}\langle (D_{1}\Phi)_{g}(\dot g(\gamma)), \kappa \rangle dV_{g}
 \neq 0. 
\end{equation}
Thus, the boundary variations $\dot g(\gamma)$ satisfying (5.30) for
 some $\kappa$ 
correspond to the cokernel. To prove (5.30), let $B(t) = \{x \in M:
 t(x) \geq t\}$ 
and $S(t) = \partial B(t) = \{x \in M: t(x) = t\}$. Apply the
 divergence theorem to 
the integral (5.30) over $B(t)$; twice for the Laplace 
term in (5.22) and once for the $\delta^{*}$ term in (5.22). Since 
$$\kappa \in Ker L \ {\rm and} \  \delta \kappa = 0,$$
it follows that the integral (5.30) reduces to an integral over the
 boundary, 
and gives
\begin{equation} \label{e5.31}
\int_{B(t)}\langle (D_{1}\Phi)_{g}(\dot g(\gamma), \kappa \rangle
 dV_{g} = 
{\tfrac{1}{2}}\int_{S(t)}(\langle \dot g(\gamma), \nabla_{N}\kappa
 \rangle  -  
\langle \nabla_{N}\dot g(\gamma), \kappa \rangle  - 
2\langle \beta(\dot g(\gamma)), \kappa(N) \rangle )dV_{S(t)}. 
\end{equation}
Of course $dV_{S(t)} = t^{-n}dV_{\gamma} + O(t^{-(n-1)})$. By (5.25)
 the last 
term in (5.31) is then $O(t)$ and so may be ignored. Let
\begin{equation}\label{e5.32}
\widetilde \kappa = t^{-n}\kappa,
\end{equation}
so that by (5.25), $|\widetilde \kappa|_{g}|_{S(t)} \leq C$. Setting
 $\hat \kappa 
= t^{2}\widetilde \kappa$, one has $|\hat \kappa|_{\bar g} =
 |\widetilde \kappa|_{g}$, 
and so the same is true for $|\hat \kappa|_{\bar g}$. From the
 definition (5.16), a 
straightforward computation shows that near $\partial M$, 
$$\dot g(\gamma) = t^{-2}\dot \gamma, \ \ {\rm and} \ \ \nabla_{N}\dot
 g(\gamma) = 0.$$ 
Note that $|\dot g(\gamma)|_{g} \sim 1$ as $t \rightarrow 0$. Hence, 
$$(\langle \dot g(\gamma), \nabla_{N}\kappa \rangle_{g}  -  
\langle \nabla_{N}\dot g(\gamma), \kappa \rangle)_{g}dV_{S(t)} = 
t^{2}\langle \nabla_{N}\kappa , \dot \gamma \rangle_{\gamma}dV_{S(t)} 
+ O(t)$$
$$= \langle \nabla_{N}\hat \kappa - (n-2)\hat \kappa, \dot \gamma
 \rangle_{\gamma} 
dV_{\gamma} + O(t).$$
Thus, 
\begin{equation}\label{e5.33}
\int_{B(t)}\langle (D_{1}\Phi)_{g}(\dot g(\gamma), \kappa \rangle
 dV_{g} = 
{\tfrac{1}{2}}\int_{S(t)}\langle \nabla_{N}\hat \kappa - (n-2)\hat
 \kappa, 
\dot \gamma \rangle_{\gamma}dV_{\gamma} + O(t).
\end{equation}

Now suppose, (contrary to (5.30)), 
\begin{equation}\label{e5.34}
\nabla_{N}\hat \kappa - (n-2)\hat \kappa = O(t),
\end{equation}
as forms on $(S(t), \bar g)$; note however that $\nabla$ is taken with 
respect to $g$ in (5.34). It follows from the smooth polyhomogeneity of 
$\hat \kappa$ near $\partial M$ and elementary integration that (5.34) 
gives
\begin{equation}\label{e5.35}
\kappa = o(t^{n}).
\end{equation}
The form $\kappa$ is an infinitesimal Einstein deformation, divergence-free 
by (5.26). Thus Corollary 4.4 and (5.35), together with the assumption in 
Theorem 1.3 that $\pi_{1}(M, \partial M) = 0$, imply that 
$$\kappa = 0 \ \ {\rm on} \ \ M,$$
giving a contradiction. This proves the relation (5.30). 

  The proof above shows that the form
\begin{equation}\label{e5.36}
\dot g(\gamma) = \lim_{t\rightarrow 0}\hat \kappa|_{S(t)},
\end{equation}
on $\partial M$ satisfies (5.30). The limit here exists by the 
smooth polyhomogeneity of $\kappa$ at $\partial M$. Thus, the 
space
\begin{equation}\label{e5.37}
\hat K = \{\hat \kappa = \lim_{t\rightarrow 0}t^{-(n-2)}\kappa|_{S(t)}:
\kappa \in K\},
\end{equation}
is naturally identified with the cokernel of $D\Pi_{g}$ in
 $T_{\gamma}\Met(\partial M)$. 
Note that $dim \widetilde K = dim K$ and also that the estimates (5.25)
 show that 
$\hat \kappa = \hat \kappa^{T}$ on $\partial M$. This means that 
infinitesimal deformations of the boundary metric $\gamma$ in the 
direction $\hat \kappa$, $\hat \kappa \in \hat K$, are not realized as 
$\frac{d}{ds}\Pi(g_{s})|_{s=0}$, where $g_{s}$ is a curve in $E_{AH}$ 
through $g$, i.e.~a curve of {\it global} Einstein metrics on $M$.

  On the other hand, suppose that $\kappa = \delta^{*}X$, i.e.~(5.28) holds 
for some $\kappa \in K$ and vector field $X$ on $M$ (necessarily) inducing 
a Killing field on $(\partial M, \gamma)$. Consider the local curve of metrics
\begin{equation}\label{e5.38}
g_{s} = g + s\delta^{*}(\frac{X}{t^{n}})
\end{equation}
defined in a neighborhood of $\partial M$. The curve $g_{s}$ is Einstein 
to $1^{\rm st}$ order in $s$ at $s = 0$. The induced variation of the 
boundary metric on $S(t)$ is, by construction, $(\widetilde \kappa)^{T}|_{S(t)} 
\sim \widetilde \kappa|_{S(t)}$, which, by rescaling, compactifies to 
$\hat \kappa$ at $\partial M$; here $\widetilde \kappa$ is given as in 
(5.32). Now note that the linearized divergence constraint (5.6) or (5.8) 
only involves the behavior at $\partial M$, or equivalently, the limiting 
behavior on $(S(t), \gamma_{t})$, $\gamma_{t} = \bar g|_{S(t)}$, as 
$t \rightarrow 0$. This basically shows that the constraint (5.6) may 
be solved in the direction $h_{(0)} = \hat \kappa$; a complete justification 
of this is given in the more computational proof to follow. Also, a simple 
calculation, cf.~(5.42) below, gives ${\mathcal L}_{X}\tau_{(n)} = 
{\mathcal L}_{X}g_{(n)} = 2\hat \kappa$. (The first statement follows 
since the term $r_{(n)}$ is intrinsic to the boundary metric $\gamma$, 
so that ${\mathcal L}_{X}r_{(n)} = 0$). Hence, it follows from Proposition 
5.4 that 
\begin{equation}\label{e5.39}
2\int_{\partial M}|\hat \kappa|^{2}dV_{\gamma} = \int_{\partial M}
\langle {\mathcal L}_{X}\tau_{(n)}, \hat \kappa \rangle dV_{\gamma} = 
2\int_{\partial M}\langle \delta(\tau_{(n)}'), X \rangle dV_{\gamma} = 
2\int_{\partial M}\langle \tau_{(n)}', \delta^{*}X \rangle dV_{\gamma} = 0,
\end{equation}
and thus ${\mathcal L}_{X}\tau_{(n)} = 0$ on $(\partial M, \gamma)$. 
Corollary 4.4 or Proposition 5.1 and the assumption $\pi_{1}(M, \partial M) 
= 0$ then imply that $\kappa = 0$ on $M$, so that $X$ is a Killing field on 
$M$. This completes the first proof of Theorem 1.3. {\endproof}

   From the converse part of Proposition 5.4, one also obtains:
\begin{corollary}\label{c5.5}
Let $g$ be a conformally compact Einstein metric on a compact manifold 
$M$ with $C^{\infty}$ boundary metric $\gamma$. Then the linearized 
divergence constraint equation {\rm (5.6)} is always solvable on 
$(\partial M, \gamma)$, i.e.~the map $\pi$ in {\rm (5.5)} is locally 
surjective at $(\gamma, \tau_{(n)})$. 
\end{corollary}

  It is useful and of interest to give another, direct computational 
proof of Theorem 1.3, without using the identification (5.37) as the 
cokernel of $D\Pi$. The basic idea is to compute as in Proposition 5.4 
on $(S(t), g_{t})$, with $A - Hg_{t}$ in place of $\tau_{(n)}$, and 
then pass to the limit on $\partial M$. Throughout the proof, we 
assume (5.28) holds. 

   Before starting the proof per se, we note that the estimates (5.25) 
and (5.28) imply that $X$ is tangential, i.e.~tangential to $(S(t), g)$, 
to high order, in that 
\begin{equation}\label{e5.40}
\langle X, N \rangle = O(t^{n+1+\mu}).
\end{equation}
To see this, one has $(\delta^{*}X)(N, N) = \langle \nabla_{N}X, 
N \rangle = N \langle X, N \rangle$. Thus (5.40) follows from (5.25) 
and the claim that $\langle X, N \rangle = 0$ on $\partial M$. To 
prove the latter, consider the compactified metric $\bar g = t^{2}g$. 
One has ${\mathcal L}_{X}\bar g = {\mathcal L}_{X}(t^{2}g) = 
2\frac{X(t)}{t}\bar g + O(t^{n})$. Thus for the induced metric 
$\gamma$ on $\partial M$, ${\mathcal L}_{X}\gamma = 
2\lambda \gamma$, where $\lambda = \lim_{t\rightarrow 0}\frac{X(t)}{t}$. 
Since $X$ is a Killing field on $(\partial M, \gamma)$, this gives
 $\lambda = 0$, which is equivalent to the statement that 
$\lim_{t\rightarrow 0}\langle X, N \rangle_{g} = 0$. Note also 
that since $X$ is smooth up to $\partial M$, $|X|_{g} = O(t^{-1})$. 

  We claim also that
\begin{equation}\label{e5.41}
[X, N] = O(t^{n+1}),
\end{equation}
in norm. First, $\langle [X, N], N\rangle = \langle \nabla_{X}N -
 \nabla_{N}X, N \rangle 
= - (\delta^{*}X)(N,N) = O(t^{n+1+\mu})$. On the other hand, on
 tangential $g$-unit vectors 
$Y$, $\langle [X, N], Y\rangle = \langle \nabla_{X}N - \nabla_{N}X, Y
 \rangle 
\sim \langle \nabla_{X}N, Y \rangle - 2(\delta^{*}X)(N,Y) + \langle
 \nabla_{Y}X, N \rangle 
\sim - 2(\delta^{*}X)(N,Y) =  O(t^{n+1})$, as claimed. Here $\sim$
 denotes equality modulo 
terms of order $o(t^{n})$. We have also used the fact that $\langle
 \nabla_{X}N, Y \rangle 
+ \langle \nabla_{Y}X, N \rangle \sim X\langle N, Y \rangle = 0$. 

   Now, to begin the proof itself, (assuming (5.28)), as above write
$$g_{s} = g + s\kappa + O(s^{2}) = g + s\delta^{*}X + O(s^{2}).$$
If $t_{s}$ is the geodesic defining function for $g_{s}$, (with
 boundary metric 
$\gamma$), then the Fefferman-Graham expansion gives $\bar g_{s} =
 dt_{s}^{2} + 
(\gamma + t_{s}^{2}g_{(2),s} + \dots + t_{s}^{n}g_{(n),s}) +
 O(t^{n+1})$. The 
estimate (5.40) implies that $t_{s} = t + sO(t^{n+2+\alpha}) +
 O(s^{2})$, 
so that modulo lower order terms, we may view $t_{s} \sim t$. Taking
 the derivative 
of the FG expansion with respect to $s$ at $s = 0$, and using the fact
 that 
$X$ is Killing on $(\partial M, \gamma)$, together with the fact that
 the lower 
order terms $g_{(k)}$, $k < n$, are determined by $\gamma$, it follows
 that, for 
$\hat \kappa$ as in (5.37), 
\begin{equation}\label{e5.42}
\hat \kappa = {\tfrac{1}{2}}{\mathcal L}_{X}g_{(n)},
\end{equation}
at $\partial M$. Here both $\hat \kappa$ and ${\mathcal L}_{X}g_{(n)}$
 are viewed 
as forms on $(\partial M, \gamma)$. 

   Next, we claim that on $(S(t), g_{t})$, 
\begin{equation}\label{e5.43}
{\mathcal L}_{X}A = -{\tfrac{n-2}{2}}t^{n-2}{\mathcal L}_{X}g_{(n)} +
 O(t^{n-1}),
\end{equation}
To see this, one has $A = \frac{1}{2}{\mathcal L}_{N}g = 
-\frac{1}{2}{\mathcal L}_{t\partial t}g = -\frac{1}{2}{\mathcal
 L}_{t\partial t}(t^{-2}g_{t})$. 
But ${\mathcal L}_{t\partial t}(t^{-2}g_{t}) = \sum {\mathcal
 L}_{t\partial t}(t^{-2+k}g_{(k)}) = 
\sum (k-2)t^{k-2}g_{(k)}$. The same reasoning as before then gives
 (5.43). 

  Given these results, we now compute
$$\int_{S(t)}\langle {\mathcal L}_{X}(A - Hg_{t}), \widetilde \kappa
 \rangle_{g_{t}} dV_{S(t)};$$
compare with the left side of (5.9). First, by (5.43),
$$\int_{S(t)}\langle {\mathcal L}_{X}A , \widetilde \kappa
 \rangle_{g_{t}} dV_{S(t)} = 
-{\tfrac{n-2}{2}}\int_{S(t)}\langle {\mathcal L}_{X}g_{(n)}, \hat
 \kappa \rangle_{\gamma} 
dV_{\gamma} + O(t).$$
Next, one has ${\mathcal L}_{X}(Hg_{t}) = X(H)g_{t} + H{\mathcal
 L}_{X}g_{t}$. For the 
first term, $X(H) = tr {\mathcal L}_{X}A + O(t^{n}) =
 -\frac{n-2}{2}t^{n-2}tr {\mathcal L}_{X}g_{(n)} 
+ O(t^{n})$. Since $tr g_{(n)}$ is intrinsic to $\gamma$ and $X$ is
 Killing on 
$(\partial M, \gamma)$, it follows that $X(H) = O(t^{n-1})$. Also,
 $\langle g_{t}, 
\widetilde \kappa \rangle = tr^{T}\widetilde \kappa$, where $tr^{T}$ is
 the tangential 
trace. By (5.25) and the fact that $\kappa$ is trace-free, $\langle
 g_{t}, \widetilde \kappa 
\rangle = O(t^{1+\alpha})$. Hence $X(H)\langle g_{t}, \widetilde \kappa
 \rangle dV_{S(t)} = 
O(t^{\alpha})$. Similarly, from (5.41) one computes ${\mathcal
 L}_{X}g_{t} = {\mathcal L}_{X}g 
+ O(t^{n+1}) = 2t^{n}\widetilde \kappa + O(t^{n+1})$. Since $H \sim n$,
 using (5.42) this gives
$$-\int_{S(t)}\langle {\mathcal L}_{X}(Hg_{t}), \widetilde \kappa
 \rangle dV_{S(t)} = 
-n\int_{S(t)}\langle {\mathcal L}_{X}g_{(n)}, \hat \kappa
 \rangle_{\gamma} 
dV_{\gamma} + O(t^{\alpha}).$$
Combining these computations then gives
\begin{equation}\label{e5.44}
\int_{S(t)}\langle {\mathcal L}_{X}(A - Hg_{t}), \widetilde \kappa
 \rangle_{g_{t}} dV_{S(t)} = 
-({\tfrac{n-2}{2}}+n)\int_{\partial M}\langle {\mathcal L}_{X}g_{(n)}, 
\hat \kappa \rangle_{\gamma} dV_{\gamma} + o(1).
\end{equation}

  On the other hand, one may use the method of proof of Proposition 
5.4 to compute the left side of (5.44). First since on $S(t)$, $\tau 
= A - Hg_{t}$ is divergence-free, a slight extension of the calculation 
(5.9) gives, for any vector field $Y$ tangent to $S(t)$ and variation $h$ 
of $g_{t} = g|_{S(t)}$, 
\begin{equation}\label{e5.45}
\int_{S(t)}\langle {\mathcal L}_{Y}\tau, h\rangle_{g_{t}} dV_{g_{t}} = 
-2\int_{S(t)}\langle \delta'(\tau), Y \rangle dV_{g_{t}} + 
\int_{S(t)}[\delta Y \langle \tau, h \rangle + \langle \tau, 
\delta^{*}Y \rangle tr h]dV_{g_{t}}.
\end{equation}
Now let the tangential variation $h$ be given by $h = 
(\delta^{*}\frac{X}{t^{n}})^{T}$, where $\delta^{*} = \delta_{g}^{*}$. Thus 
$h = (t^{-n}\kappa)^{T} = (\widetilde \kappa)^{T}$, for $\widetilde \kappa$ 
as in (5.32). Also, set $Y = X^{T}$. Observe that the estimate (5.40) 
implies that $X$ agrees with $X^{T}$ to high degree, in that 
$X = X^{T} + O(^{n+1+\mu})$. This has the effect that one may use 
$X$ and $X^{T}$ interchangably in the computations below. For example, 
since $\kappa$ is trace-free, $\delta_{g}X = 0$ and hence a simple 
calculation shows that $\delta_{g_{t}}Y = O(t^{n+1+\mu})$. Similarly, 
$tr^{T}h = O(t^{1+\mu})$, while $\delta^{*}Y = O(t^{n})$. In particular, 
the second term on the right in \eqref{e5.45} is $O(t)$, and hence may 
be ignored. 

   Next, the deformation $h$ above is the tangential part of a (trivial) 
infinitesimal Einstein deformation, and hence the linearized divergence 
constraint (5.6) holds along $S(t)$, in the direction $h$. Arguing then 
as in the proof of Proposition 5.4, it follows from \eqref{e5.45} and 
the fact that $X^{T} \sim X$ to high order that
\begin{equation}\label{e5.46}
\int_{S(t)}\langle {\mathcal L}_{X}(A - Hg_{t}), \widetilde \kappa
 \rangle_{g_{t}} dV_{S(t)} = 
2\int_{S(t)}\langle (A - Hg_{t})', (\delta^{*}X)^{T} \rangle
 dV_{g_{t}} + O(t).
\end{equation}

Now $A' = \frac{d}{ds}(A_{g+s\widetilde \kappa}) = \frac{1}{2}
({\mathcal L}_{N}\widetilde \kappa + {\mathcal L}_{N'}g) = 
\frac{1}{2}\nabla_{N}\widetilde \kappa + \widetilde \kappa + O(t)$. 
Similarly, $(Hg_{t})' = H'g_{t} + H(g_{t})'$. The first term here, 
when paired with $(\delta^{*}X)^{T}$ and integrated, gives $O(t)$, 
while the second term is $n\widetilde \kappa$ to leading order. 
Thus one has
\begin{equation}\label{e5.47}
2\int_{S(t)}\langle (A - Hg_{t})', (\delta^{*}X)^{T} \rangle dV_{g_{t}} = 
\int_{S(t)}\langle \nabla_{N}\widetilde \kappa - (2n-2)\widetilde \kappa, 
\kappa \rangle_{g_{t}} dV_{g_{t}} + O(t).
\end{equation}
A straightforward calculation, essentially the same as that preceding (5.33) 
shows that 
\begin{equation}\label{e5.48}
\int_{S(t)}\langle \nabla_{N}\widetilde \kappa - (2n-2)\widetilde \kappa, 
\kappa \rangle_{g_{t}} dV_{g_{t}} = 
\int_{S(t)}[{\tfrac{1}{2}}N(|\hat \kappa|^{2}) - (2n-2)|\hat \kappa|^{2}]
dV_{\gamma} + O(t),
\end{equation}
where the norms on the right are with respect to $\bar g$. The first term on 
the right in \eqref{e5.48} is $O(t)$, and comparing \eqref{e5.48} with 
\eqref{e5.44}-\eqref{e5.47} shows that
$$\hat \kappa = 0$$
on $\partial M$. Hence via Corollary 4.4, $\kappa = 0$ on $M$, as before. 
This completes the second proof of Theorem 1.3. 
{\endproof}

 {\bf Proof of Corollary 1.4.}

   Suppose $(M, g)$ is a conformally compact Einstein metric with
 boundary metric 
given by the round metric $S^{n}(1)$ on $S^{n}$. Theorem 1.3 implies
 that the 
isometry group of $(M, g)$ contains the isometry group of $S^{n}$. This
 reduces 
the Einstein equations to a simple system of ODE's, and it is easily
 seen that 
the only solution is given by the Poincar\'e metric on the ball
 $B^{n+1}$. 
{\endproof}

\begin{remark}\label{r5.6}
{\rm By means of Obata's theorem [30], Theorem 1.3 remains true for 
continuous groups of conformal isometries at conformal infinity. Thus, 
the class of the round metric on $S^{n}$ is the only conformal class 
which supports an essential conformal Killing field, i.e.~a field 
which is not Killing with respect to some conformally related metric. 
Corollary 1.4 shows that any $g \in E_{AH}$ with boundary metric 
$S^{n}(1)$ is necessarily the hyperbolic metric $g_{-1}$ on 
the ball. For $g_{-1}$, it is well-known that essential conformal 
Killing fields on $S^{n}$ extend to Killing fields on 
$({\mathbb H}^{n+1}, g_{-1})$. 

  We expect that a modification of the proof of Theorem 1.3 would give
 this result 
directly, without the use of Obata's theorem. In fact, such would
 probably give 
(yet) another proof of Obata's result. }
\end{remark}

  Corollary 5.5 shows, in the global situation, that the projection
 $\pi$ of 
the constraint manifold ${\mathcal T}$ to $\Met(\partial M)$ is always
 locally surjective. 
Hence there exists a formal solution, and an exact solution in the
 analytic case, for 
any nearby boundary metric, which is defined in a neighborhood of the
 boundary. However, 
the full boundary map $\Pi$ in (5.13) or (5.14) on global metrics is
 not locally surjective 
in general; nor is it always globally surjective. 

  The simplest example of this behavior is provided by the family of
 AdS Schwarzschild 
metrics. These are metrics on ${\mathbb R}^{2}\times S^{n-1}$ of the
 form
$$g_{m} = V^{-1}dr^{2} + Vd\theta^{2} + r^{2}g_{S^{n-1}(1)},$$
where $V = V(r) = 1 + r^{2} - \frac{2m}{r^{n-2}}$. Here $m > 0$ and $r
 \in [r_{+}, \infty]$, 
where $r_{+}$ is the largest root of the equation $V(r_{+}) = 0$. The
 locus $\{r_{+} = 0\}$ 
is a totally geodesic round $S^{n-1}$ of radius $r_{+}$. Smoothness of
 the metric at 
$\{r_{+} = 0\}$ requires that the circular parameter $\theta$ runs over
 the interval 
$[0,\beta]$, where 
$$\beta = \frac{4\pi r_{+}}{nr_{+}^{2}+(n-2)}.$$
The metrics $g_{m}$ are isometrically distinct for distinct values of
 $m$, and form 
a curve in $E_{AH}$ with conformal infinity given by the conformal
 class of the 
product metric on $S^{1}(\beta)\times S^{n-1}(1)$. As $m$ ranges over
 the interval 
$(0, \infty)$, $\beta$ has a maximum value of 
$$\beta \leq \beta_{\max} = 2\pi \sqrt{(n-2)/n}.$$
As $m \rightarrow 0$ or $m \rightarrow \infty$, $\beta \rightarrow 0$. 

   Hence, the metrics $S^{1}(L)\times S^{n-1}(1)$ are not in
 $\Pi(g_{m})$ for any 
$L > \beta_{\max}$. In fact these boundary metrics are not in $Im(\Pi)$
 generally, 
for any manifold $M^{n+1}$. For Theorem 1.3 implies that any
 conformally compact 
Einstein metric with boundary metric $S^{1}(L)\times S^{n-1}(1)$ has an
 isometry 
group containing the isometry group of $S^{1}(L)\times S^{n-1}(1)$.
 This again 
reduces the Einstein equations to a system of ODE's and it is easy to
 see, (although 
we do not give the calculations here), that any such metric is an AdS
 Schwarzschild metric.

\begin{remark}\label{r5.7}
{\rm In the context of Propositions 5.1 and 5.4, it is natural to
 consider the 
issue of whether local Killing fields of $\partial M$, (i.e.~Killing
 fields defined 
on the universal cover), extend to local Killing fields of any global
 conformally 
compact Einstein metric. Note that Proposition 5.1 and Proposition 5.4
 are both 
local results, the latter by using variations $h_{(0)}$ which are of
 compact support. 
However, the linearized constraint condition (5.8) is not invariant
 under covering 
spaces; even the splitting (5.7) is not invariant under coverings,
 since a Killing 
field on a covering space need not descend to the base space. 

  We claim that local Killing fields do not extend even locally into
 the interior in 
general. As a specific example, let $N^{n+1}$ be any complete,
 geometrically finite 
hyperbolic manifold, with conformal infinity $(\partial N, \gamma)$,
 and which has 
at least one parabolic end, i.e.~a finite volume cusp end, with cross
 sections given 
by flat tori $T^{n}$. There exist many such manifolds. The metric at
 conformal infinity 
is conformally flat, so there are many local Killing fields on
 $\partial N$. For example, 
in many examples $N$ itself is a compact hyperbolic manifold. Of course
 the local 
(conformal) isometries of $\partial N$ extend here to local isometries
 of $N$. 

  However, as shown in [15], the cusp end may be capped off by Dehn
 filling with a 
solid torus, to give infinitely many distinct conformally compact
 Einstein metrics 
with the same boundary metric $(\partial N, \gamma)$. These Dehn-filled
 Einstein metrics 
cannot inherit all the local conformal symmetries of the boundary. }
\end{remark}

\begin{remark} \label{r5.8} 
{\rm We point out that Theorem 1.3 fails for complete Ricci-flat
 metrics which are ALE 
(asymptotically locally Euclidean). The simplest counterexamples are
 the family of 
Eguchi-Hanson metrics, which have boundary metric at infinity given by
 the round metric 
on $S^{3}/{\mathbb Z}_{2}$. The symmetry group of these metrics is
 strictly smaller than 
the isometry group $Isom (S^{3}/{\mathbb Z}_{2})$ of the boundary.
 Similarly, the 
Gibbons-Hawking family of metrics with boundary metric the round metric
 on 
$S^{3}/{\mathbb Z}_{k}$ have only an $S^{1}$ isometry group, much
 smaller than 
the group $Isom (S^{3}/{\mathbb Z}_{k})$. 

  This indicates that, despite a number of proposals, some important
 features of holographic 
renormalization in the AdS context cannot carry over to the
 asymptotically flat case. }
\end{remark}

\bibliographystyle{plain}

\begin{thebibliography}{WWW}

\footnotesize



\bibitem [1]{1} M. Akbar and P. D'Eath, Classical boundary-value
 problem in Riemannian 
quantum gravity and self-dual Taub-NUT-(anti)de Sitter geometries,
 Nucl.Phys. {\bf B648}, 
(2003), 397-416, gr-qc/0202073, v2.

\bibitem [2]{2} M. Akbar, Classical boundary-value problem in
 Riemannian 
quantum gravity and self-dual Taub-Bolt-anti-de Sitter geometries,
 Nucl.Phys. {\bf B663}, 
(2003), 215-230, gr-qc/0301007.

\bibitem [3]{3} M. Anderson, Boundary regularity, uniqueness and
 non-uniqueness 
for AH Einstein metrics on 4-manifolds, Advances in Math., {\bf 179},
 (2003), 205-249. 

\bibitem [4]{4} M. Anderson, A. Katsuda, Y. Kurylev, M. Lassas and M.
 Taylor, 
Boundary regularity for the Ricci equation, geometric convergence and
 Gel'fand's 
inverse boundary problem, Inventiones Math., {\bf 158}, (2004),
 261-321, math.SP/0211376.

\bibitem [5]{5} M. Anderson, On the structure of conformally compact
 Einstein 
metrics, (preprint, Feb.~04/Dec.~05), math.DG/0402198.

\bibitem [6]{6} M. Anderson, Einstein metrics with prescribed conformal
 infinity on 
4-manifolds, (preprint, May 01/March 05), math.DG/0105243.

\bibitem [7]{7} M. Anderson, P. Chrus\'ciel and E. Delay, Non-trivial,
 static, 
geodesically complete space-times with a negative cosmological
 constant, II, 
in: AdS/CFT Correspondence: Einstein Metrics and Their Conformal
 Boundaries, 
Ed. O. Biquard, IMRA Lectures, vol. 8, Euro. Math. Soc, Z\"urich,
 (2005), 165-204, 
gr-qc/0401081.  

\bibitem [8]{8} L. Andersson and V. Moncrief, Elliptic-hyperbolic
 systems and the 
Einstein equations, Ann. Henri Poincar\'e, {\bf 4}, (2003), 1-34,
 gr-qc/0110111.

\bibitem [9]{9} L. Andersson and M. Dahl, Scalar curvature rigidity for
 asymptotically 
locally hyperbolic manifolds, Ann. Glob. Anal. Geom., {\bf 16}, (1998),
 1-27. 

\bibitem [10]{10} A. Besse, Einstein Manifolds, Ergebnisse der Math.
 Series 3:10, 
Springer Verlag, New York, (1987). 

\bibitem [11]{11} O. Biquard, M\'etriques d'Einstein asymptotiquement
 symm\'etriques, 
Ast\'erisque {\bf 265}, (2000). 

\bibitem [12]{12} O. Biquard, Continuation unique \`a partir de l'infini 
conforme pour les m\'etriques d'Einstein, (preprint), arXiv:0708.4346 v2.


\bibitem [13]{13} A. Calder\'on, Uniqueness in the Cauchy problem for
 partial 
differential equations, Amer. Jour. of Math., {\bf 80}, (1958), 16-36.

\bibitem [14]{14} A. Calder\'on, Existence and uniqueness theorems for
 systems of 
partial differential equations, Proc. Symp. Fluid Dynamics and Appl. 
Math., Gordon and Breach, New York, (1962), 147-195.

\bibitem [15]{15} G. Craig, Dehn filling and asymptotically hyperbolic
 Einstein 
metrics, Comm. Anal. Geom., {\bf 14}, (2006), 725-764, math.DG/0502491.
 

\bibitem [16]{16} P. Chru\'sciel, E. Delay, J. Lee and D. Skinner,
 Boundary 
regularity of conformally compact Einstein metrics, Jour. Diff. Geom.,
 {\bf 69}, 
(2005), 111-136, math.DG/0401386.

\bibitem [17]{17} S. de Haro, K. Skenderis and S. Solodukhin,
 Holographic reconstruction of 
spacetime and renormalization in the AdS/CFT correspondence, Comm.
 Math. Phys., {\bf 217}, 
(2001), 595-622, hep-th/0002230.

\bibitem [18]{18} C. Fefferman and C.R. Graham, Conformal invariants,
 in: 
\'Elie Cartan et les Mathematiques d'Aujourd'hui, Ast\'erisque, 1985,
 Numero hors Serie, 
Soc. Math. France, Paris, 95-116. 

\bibitem [19]{19} D. Gilbarg and N. Trudinger, Elliptic Partial
 Differential 
Equations of Second Order, Second Edition, Springer Verlag, New York,
 (1983). 

\bibitem [20]{20} C.R. Graham and J.M. Lee, Einstein metrics with
 prescribed 
conformal infinity on the ball, Advances in Math., {\bf 87}, (1991),
 186-225. 

\bibitem [21]{21} C. R. Graham, Volume and area renormalization for
 conformally compact 
Einstein metrics, Rend. Circ. Mat. Palermo (2) Suppl. {\bf 63}, (2000),
 31-42, 
math.DG/0009042.

\bibitem [22]{22} S.W.Hawking and G.F.R. Ellis, The Large Scale
 Structure of 
Space-Time, Cambridge Univ. Press, (1973). 

\bibitem [23]{23} J. L. Kazdan, Unique continuation in geometry, Comm.
 Pure Appl. Math., 
{\bf 41}, (1988), 667-681. 

\bibitem [24]{24} S. Kichenassamy, On a conjecture of Fefferman and
 Graham, Adv. in Math., 
{\bf 184}, (2004), 268-288.

\bibitem [25]{25} S. Kobayashi and K. Nomizu, Foundations of
 Differential Geometry, 
Wiley-Interscience, New York, 1963. 

\bibitem [26]{26} J. M. Lee, Fredholm operators and Einstein metrics on
 conformally 
compact manifolds, Memoirs Amer. Math. Soc., {\bf 183}, (2006), No.
 864, math.DG/0105046. 

\bibitem [27]{27} R. Mazzeo, Unique continuation at infinity and
 embedded eigenvalues 
for asymptotically hyperbolic manifolds, Amer. Jour. Math., {\bf 13},
 (1991), 25-45. 

\bibitem [28]{28} R. Mazzeo and F. Pacard, Maskit combinations of
 Poincar\'e-Einstein 
metrics, Advances in Math., {\bf 204}, (2006), 379-412. 

\bibitem [29]{29} L. Nirenberg, Lectures on Partial Differential
 Equations, CBMS 
Series, No.~17, Amer. Math. Soc., Providence, RI, (1973). 

\bibitem [30]{30} M. Obata, The conjectures on conformal
 transformations of Riemannian 
manifolds, Jour. Diff. Geom., {\bf 6}, (1972), 247-258. 

\bibitem [31]{31} J. Qing, On the rigidity of conformally compact
 Einstein manifolds, 
Int. Math. Res. Not., (2003), No. 21, 1141-1153, math.DG/0305084.

\bibitem [32]{32} R. Wald, General Relativity, Univ. of Chicago Press,
 Chicago, 
(1984). 

\bibitem [33]{33} X. Wang, The mass of asymptotically hyperbolic
 manifolds, Jour. Diff. 
Geom., {\bf 57}, (2001), 273-299.

\bibitem [34]{34} E. Zehnder, Generalized implicit function theorems, in: 
Topics in Nonlinear Functional Analysis, L. Nirenberg, Courant Institute 
of Math. Sci., NYU, (1974). 


\end{thebibliography}

\begin{center}
September, 2007
\end{center}

\noindent
{\address Department of Mathematics\\
S.U.N.Y. at Stony Brook\\
Stony Brook, NY 11794-3651, USA\\
E-mail: anderson@math.sunysb.edu}

\bigskip

\noindent
{\address Institut de math\'ematiques et de mod\'elisation de
 Montpellier\\
CNRS et Universit\'e Montpellier II\\
34095 Montpellier Cedex 5, France\\
E-mail: herzlich@math.univ-montp2.fr}

\end{document}